\documentclass{gtart}

\def\ifplaintex{\expandafter\ifx\csname documentclass\endcsname\relax}

\def\gtp{{\mathsurround=0pt\it $\cal G\mskip-2mu$eometry \&\ 
$\cal T\!\!$opology $\cal P\!$ublications}}  

\def\recd{{\small Received:\qua\receiveddate\ifx\reviseddate\relax
\else\qquad Revised:\qua\reviseddate\fi\par}} 

\def\lognumber#1{\def\thelognumber{#1}}
\def\volumenumber#1{\def\thevolumenumber{#1}}
\def\volumeyear#1{\def\thevolumeyear{#1}}
\def\papernumber#1{\def\thepapernumber{#1}}
\def\pagenumbers#1#2{\def\startpage{#1}\def\finishpage{#2}}
\def\published#1{\def\publishdate{#1}}

\def\received#1{\def\receiveddate{#1}}

\def\accepted#1{\def\accepteddate{#1}}

\def\asciiaddress#1{\def\theasciiaddress{#1}}

\long\def\asciiabstract#1{\long\def\theasciiabstract{#1}}

\let\\\par\let\thelognumber\relax\let\thevolumenumber\relax
\let\thepapernumber\relax\let\thevolumeyear\relax\let\startpage\relax
\let\finishpage\relax\let\publishdate\relax\let\receiveddate\relax
\let\reviseddate\relax\let\accepteddate\relax\let\theasciititle\relax
\let\theasciiauthors\relax\let\theasciiaddress\relax
\let\theasciiabstract\relax

\let\theasciiemail\relax

\ifplaintex
\font\logobig=cmssbx10 scaled 3836
\font\logomed=cmssbx10 scaled 2557
\else
\font\logobig=cmssbx10 scaled 4200
\font\logomed=cmssbx10 scaled 2800
\fi

\long\def\makeagttitle{   
\count0=\startpage
\agt\hfill      
\hbox to 45truept{\vbox to 0pt{\vglue -13truept{\logomed A\kern -.37em{\logobig 
T}\kern -.38em G}\vss}\hss}
\break
{\small Volume \thevolumenumber\ (\thevolumeyear)
\startpage--\finishpage\nl
Published: \publishdate}

\vglue .25truein

{\parskip=0pt\leftskip 0pt plus
1fil\def\\{\par\smallskip}{\Large\bf\thetitle}\par\medskip} \vglue
0.05truein

{\parskip=0pt\leftskip 0pt plus 1fil\def\\{\par}{\sc\theauthors}
\par\medskip}%
 
\vglue 0.03truein 

{\small\leftskip 25truept\rightskip 25truept{\bf Abstract}\stdspace\theabstract

{\bf AMS Classification}\stdspace\theprimaryclass
\ifx\thesecondaryclass\relax\else; \thesecondaryclass\fi\par
{\bf Keywords}\stdspace \thekeywords\par}\vglue 7truept

}   

\ifplaintex
\font\phead=cmsl9 scaled 950
\font\pnum=cmbx10 scaled 913
\font\pfoot=cmsl9 scaled 950
\headline{\vbox to 0pt{\vskip -4.5mm\line{\small\phead\ifnum
\count0=\startpage ISSN 1472-2739 (on-line) 1472-2747 (printed)
\hfill {\pnum\folio}\else\ifodd\count0\def\\{ }%
\ifx\theshorttitle\relax\thetitle\else\theshorttitle\fi\hfill{\pnum\folio}
\else\def\\{ and }{\pnum\folio}\hfill\ifx\theshortauthors\relax\theauthors
\else\theshortauthors\fi\fi\fi}\vss}}
\footline{\vbox to 0pt{\vglue 0mm\line{\small\pfoot\ifnum\count0=\startpage
\copyright\ \gtp\hfill\else
\agt, Volume \thevolumenumber\ (\thevolumeyear)\hfill\fi}\vss}}
\else
\font=cmsl9 scaled 1050
\font\lnum=cmbx10 
\font=cmsl9 scaled 1050
\makeatletter
\def\@oddhead{{\small\lhead\ifnum\count0=\startpage ISSN 1472-2739 
(on-line) 1472-2747 (printed)\hfill {\lnum\number\count0}\else\ifodd\count0
\def\\{ }\ifx\theshorttitle\relax \thetitle \else\theshorttitle\fi\hfill
{\lnum\number\count0}\else\def\\{ and }{\lnum\number\count0}
\hfill\ifx\theshortauthors\relax 
\theauthors\else\theshortauthors\fi\fi\fi}}\def\@evenhead{\@oddhead}
\def\@oddfoot{\small\lfoot\ifnum\count0=\startpage\copyright\ \gtp\hfill\else
\agt, Volume \thevolumenumber\ (\thevolumeyear)\hfill\fi}
\def\@evenfoot{\@oddfoot}
\makeatother
\fi
\let\maketitlepage\makeagttitle
\let\makeshorttitle\maketitlepage
\let\maketitle\maketitlepage

\newwrite\gtoutfile
\long\gdef\makeheadfile{  
{\def\\{, }\def\s{ }
\immediate\openout\gtoutfile head.xxx
\immediate\write\gtoutfile{To: math@arxiv.org}
\immediate\write\gtoutfile{Subject: put OR rep NNNNN:ppppp}
\immediate\write\gtoutfile{--text follows this line--}
\immediate\write\gtoutfile{Proxy-for: \ifx\theasciiauthors\relax
\theauthors\else\theasciiauthors\fi\s<\ifx\theasciiemail\relax\theemail\else\theasciiemail\fi>}
\immediate\write\gtoutfile{\noexpand\\}
\immediate\write\gtoutfile{Authors: \ifx\theasciiauthors\relax
\theauthors\else\theasciiauthors\fi}
{\def\\{ }\immediate\write\gtoutfile{Title: \ifx\theasciititle\relax
\thetitle\else\theasciititle\fi}}
\immediate\write\gtoutfile{Subj-class: GT or SG, GR etc}
\immediate\write\gtoutfile{MSC-class: \theprimaryclass\ifx\thesecondaryclass\relax\else, \thesecondaryclass\fi}
\immediate\write\gtoutfile{Journal-ref: Algebr. Geom. Topol. \thevolumenumber\s
(\thevolumeyear) \startpage-\finishpage}
\immediate\write\gtoutfile{Comments: Published by Algebraic and
Geometric Topology at}
\immediate\write\gtoutfile{\s\s\s  http://www.maths.warwick.ac.uk/agt/AGTVol\thevolumenumber/agt-\thevolumenumber-\thepapernumber.abs.html}
\immediate\write\gtoutfile{\noexpand\\}
\immediate\write\gtoutfile{}
\ifx\theasciiabstract\relax
\immediate\write\gtoutfile{\theabstract}\else
\immediate\write\gtoutfile{\theasciiabstract}\fi
\immediate\write\gtoutfile{}
\immediate\write\gtoutfile{\noexpand\\}
\immediate\write\gtoutfile{}
\immediate\closeout\gtoutfile}}  

\def\maketitlepage{\makeagttitle\makeheadfile}
\let\makeshorttitle\maketitlepage
\let\maketitle\maketitlepage

\def\ifplaintex{\expandafter\ifx\csname documentclass\endcsname\relax}

\def\gtp{{\mathsurround=0pt\it $\cal G\mskip-2mu$eometry \&\ 
$\cal T\!\!$opology $\cal P\!$ublications}}  

\def\recd{{\small Received:\qua\receiveddate\ifx\reviseddate\relax
\else\qquad Revised:\qua\reviseddate\fi\par}} 

\def\lognumber#1{\def\thelognumber{#1}}
\def\volumenumber#1{\def\thevolumenumber{#1}}
\def\volumeyear#1{\def\thevolumeyear{#1}}
\def\papernumber#1{\def\thepapernumber{#1}}
\def\pagenumbers#1#2{\def\startpage{#1}\def\finishpage{#2}}
\def\published#1{\def\publishdate{#1}}

\def\received#1{\def\receiveddate{#1}}

\def\accepted#1{\def\accepteddate{#1}}

\def\asciiaddress#1{\def\theasciiaddress{#1}}

\long\def\asciiabstract#1{\long\def\theasciiabstract{#1}}

\let\\\par\let\thelognumber\relax\let\thevolumenumber\relax
\let\thepapernumber\relax\let\thevolumeyear\relax\let\startpage\relax
\let\finishpage\relax\let\publishdate\relax\let\receiveddate\relax
\let\reviseddate\relax\let\accepteddate\relax\let\theasciititle\relax
\let\theasciiauthors\relax\let\theasciiaddress\relax
\let\theasciiabstract\relax

\let\theasciiemail\relax

\ifplaintex
\font\logobig=cmssbx10 scaled 3836
\font\logomed=cmssbx10 scaled 2557
\else
\font\logobig=cmssbx10 scaled 4200
\font\logomed=cmssbx10 scaled 2800
\fi

\long\def\makeagttitle{   
\count0=\startpage
\agt\hfill      
\hbox to 45truept{\vbox to 0pt{\vglue -13truept{\logomed A\kern -.37em{\logobig 
T}\kern -.38em G}\vss}\hss}
\break
{\small Volume \thevolumenumber\ (\thevolumeyear)
\startpage--\finishpage\nl
Published: \publishdate}

\vglue .25truein

{\parskip=0pt\leftskip 0pt plus
1fil\def\\{\par\smallskip}{\Large\bf\thetitle}\par\medskip} \vglue
0.05truein

{\parskip=0pt\leftskip 0pt plus 1fil\def\\{\par}{\sc\theauthors}
\par\medskip}%
 
\vglue 0.03truein 

{\small\leftskip 25truept\rightskip 25truept{\bf Abstract}\stdspace\theabstract

{\bf AMS Classification}\stdspace\theprimaryclass
\ifx\thesecondaryclass\relax\else; \thesecondaryclass\fi\par
{\bf Keywords}\stdspace \thekeywords\par}\vglue 7truept

}   

\ifplaintex
\font\phead=cmsl9 scaled 950
\font\pnum=cmbx10 scaled 913
\font\pfoot=cmsl9 scaled 950
\headline{\vbox to 0pt{\vskip -4.5mm\line{\small\phead\ifnum
\count0=\startpage ISSN 1472-2739 (on-line) 1472-2747 (printed)
\hfill {\pnum\folio}\else\ifodd\count0\def\\{ }%
\ifx\theshorttitle\relax\thetitle\else\theshorttitle\fi\hfill{\pnum\folio}
\else\def\\{ and }{\pnum\folio}\hfill\ifx\theshortauthors\relax\theauthors
\else\theshortauthors\fi\fi\fi}\vss}}
\footline{\vbox to 0pt{\vglue 0mm\line{\small\pfoot\ifnum\count0=\startpage
\copyright\ \gtp\hfill\else
\agt, Volume \thevolumenumber\ (\thevolumeyear)\hfill\fi}\vss}}
\else
\font=cmsl9 scaled 1050
\font\lnum=cmbx10 
\font=cmsl9 scaled 1050
\makeatletter
\def\@oddhead{{\small\lhead\ifnum\count0=\startpage ISSN 1472-2739 
(on-line) 1472-2747 (printed)\hfill {\lnum\number\count0}\else\ifodd\count0
\def\\{ }\ifx\theshorttitle\relax \thetitle \else\theshorttitle\fi\hfill
{\lnum\number\count0}\else\def\\{ and }{\lnum\number\count0}
\hfill\ifx\theshortauthors\relax 
\theauthors\else\theshortauthors\fi\fi\fi}}\def\@evenhead{\@oddhead}
\def\@oddfoot{\small\lfoot\ifnum\count0=\startpage\copyright\ \gtp\hfill\else
\agt, Volume \thevolumenumber\ (\thevolumeyear)\hfill\fi}
\def\@evenfoot{\@oddfoot}
\makeatother
\fi
\let\maketitlepage\makeagttitle
\let\makeshorttitle\maketitlepage
\let\maketitle\maketitlepage

\newwrite\gtoutfile
\long\gdef\makeheadfile{  
{\def\\{, }\def\s{ }
\immediate\openout\gtoutfile head.xxx
\immediate\write\gtoutfile{To: math@arxiv.org}
\immediate\write\gtoutfile{Subject: put OR rep NNNNN:ppppp}
\immediate\write\gtoutfile{--text follows this line--}
\immediate\write\gtoutfile{Proxy-for: \ifx\theasciiauthors\relax
\theauthors\else\theasciiauthors\fi\s<\ifx\theasciiemail\relax\theemail\else\theasciiemail\fi>}
\immediate\write\gtoutfile{\noexpand\\}
\immediate\write\gtoutfile{Authors: \ifx\theasciiauthors\relax
\theauthors\else\theasciiauthors\fi}
{\def\\{ }\immediate\write\gtoutfile{Title: \ifx\theasciititle\relax
\thetitle\else\theasciititle\fi}}
\immediate\write\gtoutfile{Subj-class: GT or SG, GR etc}
\immediate\write\gtoutfile{MSC-class: \theprimaryclass\ifx\thesecondaryclass\relax\else, \thesecondaryclass\fi}
\immediate\write\gtoutfile{Journal-ref: Algebr. Geom. Topol. \thevolumenumber\s
(\thevolumeyear) \startpage-\finishpage}
\immediate\write\gtoutfile{Comments: Published by Algebraic and
Geometric Topology at}
\immediate\write\gtoutfile{\s\s\s  http://www.maths.warwick.ac.uk/agt/AGTVol\thevolumenumber/agt-\thevolumenumber-\thepapernumber.abs.html}
\immediate\write\gtoutfile{\noexpand\\}
\immediate\write\gtoutfile{}
\ifx\theasciiabstract\relax
\immediate\write\gtoutfile{\theabstract}\else
\immediate\write\gtoutfile{\theasciiabstract}\fi
\immediate\write\gtoutfile{}
\immediate\write\gtoutfile{\noexpand\\}
\immediate\write\gtoutfile{}
\immediate\closeout\gtoutfile}}  

\def\maketitlepage{\makeagttitle\makeheadfile}
\let\makeshorttitle\maketitlepage
\let\maketitle\maketitlepage

\def\ifplaintex{\expandafter\ifx\csname documentclass\endcsname\relax}

\def\gtp{{\mathsurround=0pt\it $\cal G\mskip-2mu$eometry \&\ 
$\cal T\!\!$opology $\cal P\!$ublications}}  

\def\recd{{\small Received:\qua\receiveddate\ifx\reviseddate\relax
\else\qquad Revised:\qua\reviseddate\fi\par}} 

\def\lognumber#1{\def\thelognumber{#1}}
\def\volumenumber#1{\def\thevolumenumber{#1}}
\def\volumeyear#1{\def\thevolumeyear{#1}}
\def\papernumber#1{\def\thepapernumber{#1}}
\def\pagenumbers#1#2{\def\startpage{#1}\def\finishpage{#2}}
\def\published#1{\def\publishdate{#1}}

\def\received#1{\def\receiveddate{#1}}

\def\accepted#1{\def\accepteddate{#1}}

\def\asciiaddress#1{\def\theasciiaddress{#1}}

\long\def\asciiabstract#1{\long\def\theasciiabstract{#1}}

\let\\\par\let\thelognumber\relax\let\thevolumenumber\relax
\let\thepapernumber\relax\let\thevolumeyear\relax\let\startpage\relax
\let\finishpage\relax\let\publishdate\relax\let\receiveddate\relax
\let\reviseddate\relax\let\accepteddate\relax\let\theasciititle\relax
\let\theasciiauthors\relax\let\theasciiaddress\relax
\let\theasciiabstract\relax

\let\theasciiemail\relax

\ifplaintex
\font\logobig=cmssbx10 scaled 3836
\font\logomed=cmssbx10 scaled 2557
\else
\font\logobig=cmssbx10 scaled 4200
\font\logomed=cmssbx10 scaled 2800
\fi

\long\def\makeagttitle{   
\count0=\startpage
\agt\hfill      
\hbox to 45truept{\vbox to 0pt{\vglue -13truept{\logomed A\kern -.37em{\logobig 
T}\kern -.38em G}\vss}\hss}
\break
{\small Volume \thevolumenumber\ (\thevolumeyear)
\startpage--\finishpage\nl
Published: \publishdate}

\vglue .25truein

{\parskip=0pt\leftskip 0pt plus
1fil\def\\{\par\smallskip}{\Large\bf\thetitle}\par\medskip} \vglue
0.05truein

{\parskip=0pt\leftskip 0pt plus 1fil\def\\{\par}{\sc\theauthors}
\par\medskip}%
 
\vglue 0.03truein 

{\small\leftskip 25truept\rightskip 25truept{\bf Abstract}\stdspace\theabstract

{\bf AMS Classification}\stdspace\theprimaryclass
\ifx\thesecondaryclass\relax\else; \thesecondaryclass\fi\par
{\bf Keywords}\stdspace \thekeywords\par}\vglue 7truept

}   

\ifplaintex
\font\phead=cmsl9 scaled 950
\font\pnum=cmbx10 scaled 913
\font\pfoot=cmsl9 scaled 950
\headline{\vbox to 0pt{\vskip -4.5mm\line{\small\phead\ifnum
\count0=\startpage ISSN 1472-2739 (on-line) 1472-2747 (printed)
\hfill {\pnum\folio}\else\ifodd\count0\def\\{ }%
\ifx\theshorttitle\relax\thetitle\else\theshorttitle\fi\hfill{\pnum\folio}
\else\def\\{ and }{\pnum\folio}\hfill\ifx\theshortauthors\relax\theauthors
\else\theshortauthors\fi\fi\fi}\vss}}
\footline{\vbox to 0pt{\vglue 0mm\line{\small\pfoot\ifnum\count0=\startpage
\copyright\ \gtp\hfill\else
\agt, Volume \thevolumenumber\ (\thevolumeyear)\hfill\fi}\vss}}
\else
\font=cmsl9 scaled 1050
\font\lnum=cmbx10 
\font=cmsl9 scaled 1050
\makeatletter
\def\@oddhead{{\small\lhead\ifnum\count0=\startpage ISSN 1472-2739 
(on-line) 1472-2747 (printed)\hfill {\lnum\number\count0}\else\ifodd\count0
\def\\{ }\ifx\theshorttitle\relax \thetitle \else\theshorttitle\fi\hfill
{\lnum\number\count0}\else\def\\{ and }{\lnum\number\count0}
\hfill\ifx\theshortauthors\relax 
\theauthors\else\theshortauthors\fi\fi\fi}}\def\@evenhead{\@oddhead}
\def\@oddfoot{\small\lfoot\ifnum\count0=\startpage\copyright\ \gtp\hfill\else
\agt, Volume \thevolumenumber\ (\thevolumeyear)\hfill\fi}
\def\@evenfoot{\@oddfoot}
\makeatother
\fi
\let\maketitlepage\makeagttitle
\let\makeshorttitle\maketitlepage
\let\maketitle\maketitlepage

\newwrite\gtoutfile
\long\gdef\makeheadfile{  
{\def\\{, }\def\s{ }
\immediate\openout\gtoutfile head.xxx
\immediate\write\gtoutfile{To: math@arxiv.org}
\immediate\write\gtoutfile{Subject: put OR rep NNNNN:ppppp}
\immediate\write\gtoutfile{--text follows this line--}
\immediate\write\gtoutfile{Proxy-for: \ifx\theasciiauthors\relax
\theauthors\else\theasciiauthors\fi\s<\ifx\theasciiemail\relax\theemail\else\theasciiemail\fi>}
\immediate\write\gtoutfile{\noexpand\\}
\immediate\write\gtoutfile{Authors: \ifx\theasciiauthors\relax
\theauthors\else\theasciiauthors\fi}
{\def\\{ }\immediate\write\gtoutfile{Title: \ifx\theasciititle\relax
\thetitle\else\theasciititle\fi}}
\immediate\write\gtoutfile{Subj-class: GT or SG, GR etc}
\immediate\write\gtoutfile{MSC-class: \theprimaryclass\ifx\thesecondaryclass\relax\else, \thesecondaryclass\fi}
\immediate\write\gtoutfile{Journal-ref: Algebr. Geom. Topol. \thevolumenumber\s
(\thevolumeyear) \startpage-\finishpage}
\immediate\write\gtoutfile{Comments: Published by Algebraic and
Geometric Topology at}
\immediate\write\gtoutfile{\s\s\s  http://www.maths.warwick.ac.uk/agt/AGTVol\thevolumenumber/agt-\thevolumenumber-\thepapernumber.abs.html}
\immediate\write\gtoutfile{\noexpand\\}
\immediate\write\gtoutfile{}
\ifx\theasciiabstract\relax
\immediate\write\gtoutfile{\theabstract}\else
\immediate\write\gtoutfile{\theasciiabstract}\fi
\immediate\write\gtoutfile{}
\immediate\write\gtoutfile{\noexpand\\}
\immediate\write\gtoutfile{}
\immediate\closeout\gtoutfile}}  

\def\maketitlepage{\makeagttitle\makeheadfile}
\let\makeshorttitle\maketitlepage
\let\maketitle\maketitlepage

\def\ifplaintex{\expandafter\ifx\csname documentclass\endcsname\relax}

\def\gtp{{\mathsurround=0pt\it $\cal G\mskip-2mu$eometry \&\ 
$\cal T\!\!$opology $\cal P\!$ublications}}  

\def\recd{{\small Received:\qua\receiveddate\ifx\reviseddate\relax
\else\qquad Revised:\qua\reviseddate\fi\par}} 

\def\lognumber#1{\def\thelognumber{#1}}
\def\volumenumber#1{\def\thevolumenumber{#1}}
\def\volumeyear#1{\def\thevolumeyear{#1}}
\def\papernumber#1{\def\thepapernumber{#1}}
\def\pagenumbers#1#2{\def\startpage{#1}\def\finishpage{#2}}
\def\published#1{\def\publishdate{#1}}

\def\received#1{\def\receiveddate{#1}}

\def\accepted#1{\def\accepteddate{#1}}

\def\asciiaddress#1{\def\theasciiaddress{#1}}

\long\def\asciiabstract#1{\long\def\theasciiabstract{#1}}

\let\\\par\let\thelognumber\relax\let\thevolumenumber\relax
\let\thepapernumber\relax\let\thevolumeyear\relax\let\startpage\relax
\let\finishpage\relax\let\publishdate\relax\let\receiveddate\relax
\let\reviseddate\relax\let\accepteddate\relax\let\theasciititle\relax
\let\theasciiauthors\relax\let\theasciiaddress\relax
\let\theasciiabstract\relax

\let\theasciiemail\relax

\ifplaintex
\font\logobig=cmssbx10 scaled 3836
\font\logomed=cmssbx10 scaled 2557
\else
\font\logobig=cmssbx10 scaled 4200
\font\logomed=cmssbx10 scaled 2800
\fi

\long\def\makeagttitle{   
\count0=\startpage
\agt\hfill      
\hbox to 45truept{\vbox to 0pt{\vglue -13truept{\logomed A\kern -.37em{\logobig 
T}\kern -.38em G}\vss}\hss}
\break
{\small Volume \thevolumenumber\ (\thevolumeyear)
\startpage--\finishpage\nl
Published: \publishdate}

\vglue .25truein

{\parskip=0pt\leftskip 0pt plus
1fil\def\\{\par\smallskip}{\Large\bf\thetitle}\par\medskip} \vglue
0.05truein

{\parskip=0pt\leftskip 0pt plus 1fil\def\\{\par}{\sc\theauthors}
\par\medskip}%
 
\vglue 0.03truein 

{\small\leftskip 25truept\rightskip 25truept{\bf Abstract}\stdspace\theabstract

{\bf AMS Classification}\stdspace\theprimaryclass
\ifx\thesecondaryclass\relax\else; \thesecondaryclass\fi\par
{\bf Keywords}\stdspace \thekeywords\par}\vglue 7truept

}   

\ifplaintex
\font\phead=cmsl9 scaled 950
\font\pnum=cmbx10 scaled 913
\font\pfoot=cmsl9 scaled 950
\headline{\vbox to 0pt{\vskip -4.5mm\line{\small\phead\ifnum
\count0=\startpage ISSN 1472-2739 (on-line) 1472-2747 (printed)
\hfill {\pnum\folio}\else\ifodd\count0\def\\{ }%
\ifx\theshorttitle\relax\thetitle\else\theshorttitle\fi\hfill{\pnum\folio}
\else\def\\{ and }{\pnum\folio}\hfill\ifx\theshortauthors\relax\theauthors
\else\theshortauthors\fi\fi\fi}\vss}}
\footline{\vbox to 0pt{\vglue 0mm\line{\small\pfoot\ifnum\count0=\startpage
\copyright\ \gtp\hfill\else
\agt, Volume \thevolumenumber\ (\thevolumeyear)\hfill\fi}\vss}}
\else
\font=cmsl9 scaled 1050
\font\lnum=cmbx10 
\font=cmsl9 scaled 1050
\makeatletter
\def\@oddhead{{\small\lhead\ifnum\count0=\startpage ISSN 1472-2739 
(on-line) 1472-2747 (printed)\hfill {\lnum\number\count0}\else\ifodd\count0
\def\\{ }\ifx\theshorttitle\relax \thetitle \else\theshorttitle\fi\hfill
{\lnum\number\count0}\else\def\\{ and }{\lnum\number\count0}
\hfill\ifx\theshortauthors\relax 
\theauthors\else\theshortauthors\fi\fi\fi}}\def\@evenhead{\@oddhead}
\def\@oddfoot{\small\lfoot\ifnum\count0=\startpage\copyright\ \gtp\hfill\else
\agt, Volume \thevolumenumber\ (\thevolumeyear)\hfill\fi}
\def\@evenfoot{\@oddfoot}
\makeatother
\fi
\let\maketitlepage\makeagttitle
\let\makeshorttitle\maketitlepage
\let\maketitle\maketitlepage

\newwrite\gtoutfile
\long\gdef\makeheadfile{  
{\def\\{, }\def\s{ }
\immediate\openout\gtoutfile head.xxx
\immediate\write\gtoutfile{To: math@arxiv.org}
\immediate\write\gtoutfile{Subject: put OR rep NNNNN:ppppp}
\immediate\write\gtoutfile{--text follows this line--}
\immediate\write\gtoutfile{Proxy-for: \ifx\theasciiauthors\relax
\theauthors\else\theasciiauthors\fi\s<\ifx\theasciiemail\relax\theemail\else\theasciiemail\fi>}
\immediate\write\gtoutfile{\noexpand\\}
\immediate\write\gtoutfile{Authors: \ifx\theasciiauthors\relax
\theauthors\else\theasciiauthors\fi}
{\def\\{ }\immediate\write\gtoutfile{Title: \ifx\theasciititle\relax
\thetitle\else\theasciititle\fi}}
\immediate\write\gtoutfile{Subj-class: GT or SG, GR etc}
\immediate\write\gtoutfile{MSC-class: \theprimaryclass\ifx\thesecondaryclass\relax\else, \thesecondaryclass\fi}
\immediate\write\gtoutfile{Journal-ref: Algebr. Geom. Topol. \thevolumenumber\s
(\thevolumeyear) \startpage-\finishpage}
\immediate\write\gtoutfile{Comments: Published by Algebraic and
Geometric Topology at}
\immediate\write\gtoutfile{\s\s\s  http://www.maths.warwick.ac.uk/agt/AGTVol\thevolumenumber/agt-\thevolumenumber-\thepapernumber.abs.html}
\immediate\write\gtoutfile{\noexpand\\}
\immediate\write\gtoutfile{}
\ifx\theasciiabstract\relax
\immediate\write\gtoutfile{\theabstract}\else
\immediate\write\gtoutfile{\theasciiabstract}\fi
\immediate\write\gtoutfile{}
\immediate\write\gtoutfile{\noexpand\\}
\immediate\write\gtoutfile{}
\immediate\closeout\gtoutfile}}  

\def\maketitlepage{\makeagttitle\makeheadfile}
\let\makeshorttitle\maketitlepage
\let\maketitle\maketitlepage

\lognumber{28}
\volumenumber{2}
\volumeyear{2002}
\papernumber{28}
\published{25 July 2002}
\pagenumbers{591}{647}
\received{29 January 2002}
\accepted{25 June 2002}

\usepackage{amsmath,amssymb}
\usepackage[all]{xy}

\newtheorem{thm}{Theorem}[section]
\newtheorem{lem}[thm]{Lemma}
\newtheorem{cor}[thm]{Corollary}
\newtheorem{prop}[thm]{Proposition}

\theoremstyle{definition}

\newtheorem{defn}[thm]{Definition}
\newtheorem{defns}[thm]{Definitions}

\newtheorem{notation}[thm]{Notation}
\newtheorem{ex}[thm]{Example}

\theoremstyle{remark}

\newtheorem{rem}[thm]{Remark}

\numberwithin{equation}{section}

\newcommand{\thmref}[1]{Theorem~\ref{#1}}
\newcommand{\corref}[1]{Corollary~\ref{#1}}
\newcommand{\secref}[1]{\S\ref{#1}}

\newcommand{\propref}[1]{Proposition~\ref{#1}}
\newcommand{\lemref}[1]{Lemma~\ref{#1}}
\newcommand{\appref}[1]{Appendix~\ref{#1}}
\newcommand{\exref}[1]{Example~\ref{#1}}

\newcommand{\hocolim}{\operatorname*{hocolim}}
\newcommand{\holim}{\operatorname*{holim}}
\newcommand{\colim}{\operatorname*{colim}}
\newcommand{\Hom}{\operatorname{Hom}}
\newcommand{\Map}{\operatorname{Map}}
\newcommand{\MapS}{\operatorname{Map_{\mathcal S}}}
\newcommand{\MapT}{\operatorname{Map_{\mathcal T}}}

\newcommand{\C}{{\mathcal  C}}
\newcommand{\Sp}{{\mathcal S}}
\newcommand{\T}{{\mathcal  T}}
\newcommand{\J}{{\mathcal  J}}
\newcommand{\E}{{\mathcal  E}}
\newcommand{\U}{{\mathcal  U}}
\newcommand{\PP}{{\mathcal  P}}

\newcommand{\add}{{\mathbf  a}}

\newcommand{\Z}{{\mathbb  Z}}
\newcommand{\N}{{\mathbb  N}}

\newcommand{\Sinfty}{\Sigma^{\infty}}

\newcommand{\sm}{\wedge}

\newcommand{\da}{\downarrow}

\newcommand{\ra}{\rightarrow}
\newcommand{\xra}{\xrightarrow}
\newcommand{\la}{\leftarrow}
\newcommand{\xla}{\xleftarrow}

\newcommand{\hra}{\hookrightarrow}

\begin{document}

\title[Structure in resolutions of mapping spaces]{Product and other fine structure in polynomial\\resolutions of mapping spaces}

\authors{Stephen T. Ahearn\\Nicholas J. Kuhn}

\address{Department of Mathematics, Macalester College\\St.Paul, MN 55105, USA}    
\email{ahearn@macalester.edu, njk4x@virginia.edu}
\secondaddress{Department of Mathematics, University of Virginia\\Charlottesville, VA 22903, USA}    

\asciiaddress{Department of Mathematics, Macalester College\\St.Paul, MN 55105, USA\\and\\Department of Mathematics, University of Virginia\\Charlottesville, VA 22903, USA}

\begin{abstract}
Let $\operatorname{Map_{\mathcal T}}(K,X)$ denote the mapping space of continuous based functions between two based spaces $K$ and $X$.  If $K$ is a fixed finite complex, Greg Arone has recently given an explicit model for the Goodwillie tower of the functor sending a space $X$ to the suspension spectrum $\Sigma^{\infty} \operatorname{Map_{\mathcal T}}(K,X)$.  

Applying a generalized homology theory $h_*$ to this tower yields a spectral sequence, and this will converge strongly to $h_*(\operatorname{Map_{\mathcal T}}(K,X))$ under suitable conditions, e.g.~if $h_*$ is connective and $X$ is at least $\text{dim } K$ connected.  Even when the convergence is more problematic, it appears the spectral sequence can still shed considerable light on $h_*(\operatorname{Map_{\mathcal T}}(K,X))$.  Similar comments hold when a cohomology theory is applied.

In this paper we study how various important natural constructions on mapping spaces induce extra structure on the towers. This leads to useful interesting additional structure in the associated spectral sequences.  For example, the diagonal on $\operatorname{Map_{\mathcal T}}(K,X)$ induces a `diagonal' on the associated tower.  After applying any cohomology theory with products $h^*$, the resulting spectral sequence is then a spectral sequence of differential graded algebras.  The product on the $E_{\infty}$--term corresponds to the cup product in $h^*(\operatorname{Map_{\mathcal T}}(K,X))$ in the usual way, and the product on the $E_1$--term  is described in terms of group theoretic transfers.

We use explicit equivariant S--duality maps to show that, when $K$ is the sphere $S^n$, our constructions at the fiber level have descriptions in terms of the Boardman--Vogt little $n$--cubes spaces.  We are then able to identify, in a computationally useful way, the Goodwillie tower of the functor from spectra to spectra sending a spectrum $X$ to $\Sigma^{\infty} \Omega^{\infty} X$.
\end{abstract}
\asciiabstract{ 
Let Map_T(K,X) denote the mapping space of continuous based functions
between two based spaces K and X.  If K is a fixed finite complex,
Greg Arone has recently given an explicit model for the Goodwillie
tower of the functor sending a space X to the suspension spectrum
\Sigma^\infty Map_T(K,X).  

Applying a generalized homology theory h_* to this tower yields a
spectral sequence, and this will converge strongly to h_*(Map_T(K,X))
under suitable conditions, e.g. if h_* is connective and X is at least
dim K connected.  Even when the convergence is more problematic, it
appears the spectral sequence can still shed considerable light on
h_*(Map_T(K,X)).  Similar comments hold when a cohomology theory is
applied.

In this paper we study how various important natural constructions on
mapping spaces induce extra structure on the towers. This leads to
useful interesting additional structure in the associated spectral
sequences.  For example, the diagonal on Map_T(K,X) induces a
`diagonal' on the associated tower.  After applying any cohomology
theory with products h^*, the resulting spectral sequence is then a
spectral sequence of differential graded algebras.  The product on the
E_\infty-term corresponds to the cup product in h^*(Map_T(K,X)) in the
usual way, and the product on the E_1-term is described in terms of
group theoretic transfers.

We use explicit equivariant S-duality maps to show that, when K is the
sphere S^n, our constructions at the fiber level have descriptions in
terms of the Boardman-Vogt little n-cubes spaces.  We are then able to
identify, in a computationally useful way, the Goodwillie tower of the
functor from spectra to spectra sending a spectrum X to
\Sigma^\infty\Omega^\infty X.}

\maketitle

\section{Introduction} \label{introduction}

Let $\MapT(K,X)$ denote the mapping space of continuous based maps between two based spaces $K$ and $X$. To compute its homology or cohomology with respect to any generalized theory, it suffices to consider the suspension spectrum $\Sigma^{\infty}\MapT(K,X)_+$, where $Z_+$ denotes the union of a space $Z$ with a disjoint basepoint.  

If one fixes $K$ and lets $X$ vary, one gets a functor from spaces to spectra.  Assuming, as we will also do from now on, that $K$ is a finite CW complex, G.~Arone \cite{arone} has recently studied this functor from the point of T.~Goodwillie's calculus of functors \cite{goodwillie1,goodwillie2,goodwillie3}.  He defines a very explicit natural tower $P^K(X)$ of fibrations of spectra under $\Sigma^{\infty}\MapT(K,X)_+$,
\begin{equation*}
\xymatrix{
&&& \vdots \ar[d] \\
&&& P^K_2(X) \ar[d] \\
&&& P^K_1(X) \ar[d] \\
\Sinfty \MapT(K,X)_+ \ar[rrr]  \ar[urrr] \ar[uurrr] &&& P^K_0(X), 
}
\end{equation*}
and shows that the connectivity of the maps 
$$e_k^K(X): \Sigma^{\infty}\MapT(K,X)_+ \ra P^K_k(X)$$ 
increases linearly with $k$ as long as the dimension of $K$ is no more than the connectivity of $X$.  The $k^{th}$ fiber $F^K_k(X)$ of the tower is shown to be naturally weakly equivalent to a homotopy orbit spectrum:
\begin{equation} \label{fiber} F^K_k(X) \simeq \MapS(K^{(k)}, X^{\sm k})_{h \Sigma_k}.
\end{equation}
Here $K^{(k)} = K^{\sm k}/\Delta_k(K)$, the quotient of the $k$--fold smash product $K^{\sm k}$ by the fat diagonal $\Delta_k(K)$, and $\MapS(K^{(k)}, X^{\sm k})$ denotes the spectrum of stable maps from $K^{(k)}$ to $ X^{\sm k}$, a spectrum with an action of the $k^{th}$ symmetric group $\Sigma_k$.
Since this is a homogeneous polynomial functor of degree $k$, Arone has identified the Goodwillie tower of $\Sigma^{\infty}\MapT(K,X)_+$.

Applying a generalized homology theory $h_*$ to this tower yields a (left half plane) spectral sequence, and this will converge strongly to $h_*(\MapT(K,X))$ under suitable conditions, e.g.\ if $h_*$ is connective and $X$ is at least $\text{dim } K$ connected.  Even when the convergence is more problematic, it appears the spectral sequence can still shed considerable light on $h_*(\MapT(K,X))$.  Similar comments hold when a cohomology theory is applied.

When $K$ is the circle $S^1$, and the homology theory is ordinary, one can show that the resulting spectral sequence is  the classical Eilenberg--Moore spectral sequence. For other $K$, it appears that the Arone spectral sequences are organized more usefully than the older Anderson spectral sequence \cite{anderson} for computing the homology and cohomology of $\MapT(K,X)$.\footnote{We note that \cite{bg} suggests that the two spectral sequences are related.}

For the deepest applications of essentially any interesting spectral sequence, one uses additional structure that the spectral sequences carries.
It is the purpose of this paper to study various geometric properties of the towers $P^K(X)$ which lead to such interesting additional structure in their associated spectral sequences.  For example, we construct a `diagonal' on $P^K(X)$.  After applying any cohomology theory with products $h^*$, the resulting spectral sequence will then be a spectral sequence of differential graded algebras.  The product on the $E_{\infty}$--term will correspond to the cup product in $h^*(\MapT(K,X))$ in the usual way, and the product on the $E_1$--term will be described in terms of group theoretic transfers.

Perhaps the towers of greatest interest are those when $K = S^n$, the $n$--sphere.  We combine (\ref{fiber}) with an explicit unstable $\Sigma_k$--equivariant S--duality map
$$ \delta(n,k):\C(n,k)_+ \sm S^{n(k)} \ra S^{nk},$$
to construct an explicit natural weak homotopy equivalence   
\begin{equation} \label{Sn fiber1} F^{S^n}_k(X) \simeq (\C(n,k)_+ \sm \MapS(S^{nk}, X^{\sm k}))_{h\Sigma_k}.
\end{equation}
Here $\C(n,k)$ is the Boardman--Vogt space of $k$ disjoint little $n$--cubes in a big $n$--cube \cite{may}.

In terms of the extended power constructions of \cite{lmms}, this last equivalence yields a weak equivalence
\begin{equation} \label{Sn fiber2} F^{S^n}_k(X) \simeq  \C(n,k)_+ \sm_{\Sigma_k} (\Sigma^{-n} X)^{\sm k}.
\end{equation}
Here $\Sigma^{-n} X$ denotes the $n^{th}$ desuspension of the suspension spectrum of $X$.  

Using either (\ref{Sn fiber1}) or (\ref{Sn fiber2}), our general structure theorems for $P^K(X)$ simplify in  nice ways when specialized to $P^{S^n}(X)$. This leads to the spectral sequences for computing $h_*(\Omega^n X)$ having lots of extra algebraic structure that can be related to classical calculations, together with statements about how these spectral sequences are related as $n$ varies.

We will not give applications in this paper.  However, in work that will appear elsewhere, the second author has used just a small part of the structure in the spectral sequences for computing $H^*(\Omega^n X; \Z/2)$ to simplify the proof of some topological nonrealization results of L.Schwartz \cite{schwartz}.   This structure also appears to be a reflection of structure in spectral sequences for calculating versions of higher Topological Hochschild Homology (see \cite{k3}).

\subsection{The Smashing Theorem}

Our first result is our simplest and most expected.  It arises from  the natural map between function spaces
$$ \eta: \MapT(L,X) \ra \MapT(K\sm L, K\sm X)$$
that one gets by smashing with the identity map of $K$. \\

\begin{thm}   There are natural maps of towers
$$\eta:P^L(X) \ra P^{K \sm L}(K \sm X)$$
 with the following properties. 
\begin{enumerate}
\item There is a commutative diagram of spectra:
$$
\xymatrix{
\Sinfty\MapT(L,X)_+    \ar[d]^{\Sinfty \eta_+} \ar[rr]^-{e^L(X)} & &
P^L(X)   \ar[d]^{\eta}  \\
\Sinfty\MapT(K \sm L,K \sm X)_+ \ar[rr]^-{e^{K \sm L}(K \sm X)} & &
P^{K \sm L}(K\sm X). 
}
$$
\item The induced map on $k^{th}$ fibers 
$$ F_k^{L}(X) \ra F_k^{K \sm L}(K \sm X)$$
is naturally equivalent  to the composite 
\begin{equation*}
\begin{split}
\Map_{\Sp}(L^{( k)}, X^{\sm k})_{h\Sigma_k} & 
\xra{\eta} \Map_{\Sp}(K^{\sm k} \sm L^{( k)}, K^{\sm k} \sm X^{\sm k})_{h\Sigma_k} \\
& \xra{p^*} \Map_{\Sp}((K\sm L)^{( k)}, (K \sm X)^{\sm k})_{h\Sigma_k},
\end{split} 
\end{equation*}
where $p:(K \sm L)^{(k)} \ra K^{\sm k}\sm L^{(k)}$ is the $\Sigma_k$--equivariant projection.
\end{enumerate}
\end{thm}

\begin{cor}  \label{Sn smash cor}
There is a natural map of towers
$$ \eta: P^{S^m}(X) \ra P^{S^{m+n}}(\Sigma^n X) $$
under 
$$\Sinfty \eta_+: \Sinfty (\Omega^m X)_+ \ra 
\Sinfty (\Omega^{m+n} \Sigma^{n}X)_+,$$ such that the associated map on $k^{th}$ fibers is equivalent to the map 
$$\C(m,k)_+ \sm_{\Sigma_k} (\Sigma^{-m}X)^{\sm k} \ra \C(m+n,k)_+ \sm_{\Sigma_k} (\Sigma^{-m}X)^{\sm k}$$
induced by the $\Sigma_k$--equivariant inclusion $\C(m,k) \hra \C(m+n,k)$. \\
\end{cor}

We have listed this theorem and corollary first because it allows us to extend the definition of our towers for $\Sinfty \Omega^n X$, with $X$ a space, to towers for  $\Sinfty \Omega^{\infty} X$, with $X$ a spectrum.  Let the spaces $\{X_n\}$, $n \geq 0$, be the spaces in the spectrum $X$, so that $\Omega^n X_n = \Omega^{\infty} X$ for all $n$.  Then {\em define} $P^{S^{\infty}}(X)$ to be the hocolimit over $n$ of the maps of towers
$$ P^{S^n}(X_n) \ra P^{S^{n+1}}(\Sigma X_n)\ra P^{S^{n+1}}(X_{n+1})$$
where the first map is given by the theorem and the second by the spectrum structure maps.  Recalling that \ $\displaystyle \hocolim_n \Sigma^{-n} \Sinfty X_n$ \ is naturally equivalent to $X$, the maps $\Sinfty \Omega^n X_{n+} \ra P^{S^n}(X_n)$ induce maps
$$ \Sinfty \Omega^{\infty} X_+ \ra P^{S^{\infty}}(X).$$
We deduce the following. \\

\begin{cor}  The $k^{th}$ fiber of the tower $P^{S^{\infty}}(X)$ is naturally equivalent to the $k^{th}$ extended power
$$ \C(\infty,k)_+ \sm_{\Sigma_k} X^{\sm k} \simeq (X^{\sm k})_{h\Sigma_k}.$$
If $X$ is $0$--connected, then the connectivity of the maps $ \Sinfty \Omega^{\infty} X_+ \ra P^{S^{\infty}}_k(X)$ increases linearly with $k$. \\
\end{cor}

This identification of both the fibers and convergence of the Goodwillie tower of $\Sinfty \Omega^{\infty}: \text{Spectra} \ra \text{Spectra}$ has been observed previously by other people, e.g.\ Goodwillie, Arone, and R.McCarthy (see the comments at the beginning of \cite{mccarthy2}).  However, by constructing it in this way, our later structure theorems for $P^{S^n}(X)$ will immediately imply analogous results about $P^{S^{\infty}}(X)$, and thus results about spectral sequences for computing $h_*(\Omega^{\infty}X)$. \\

\subsection{The Product and Diagonal Theorems}

Our next results are consequences of our study of the map of towers associated to the natural homeomorphisms of function spaces
$$ \MapT(K\vee L, X) = \MapT(K,X) \times \MapT(L,X).$$
To state these, we need to introduce a bit of the language one would use in defining the homotopy category of functors, and also describe an appropriate sort of  completed smash product of towers of spectra.

For the former, given two functors $F$ and $G$ from pointed spaces to spectra, a {\em weak} natural transformation $h: F \rightarrow G$ will be a triple $(H,f,g)$, with $H$ a functor from spaces to spectra, $g: H \rightarrow G$ a natural transformation, and $f: H \rightarrow F$ a natural transformation such that $f(X): H(X) \rightarrow F(X)$ is a weak homotopy equivalence for all $X$. If $g(X)$ is also a weak homotopy equivalence for all $X$, then we say that $h$ is a weak natural equivalence. Note that, if $F$ and $G$ are homotopy functors, then a weak natural transformation $h: F \rightarrow G$ induces a well defined natural transformation in the homotopy category: $h(X) = g(X) \circ f(X)^{-1} \in [F(X), G(X)]$.  Furthermore, using homotopy pullbacks, one can define the composition of weak natural transformations.

Now we need to define the smash product of two towers of spectra.  If $P$ and $Q$ are two towers of spectra, let $P\sm Q$ be the tower with 
$$ (P \sm Q)_k = \holim_{i+j \leq k} P_i \sm Q_j.$$
Let $F_k(P)$ denote the homotopy fiber of $P_k \ra P_{k-1}$.  As will be noted in \secref{towers of spectra}, there is a weak natural equivalence
$$ F_k(P\sm Q) \simeq \prod_{i+j = k} F_i(P) \sm F_j(Q).$$

\begin{thm}   There are natural weak homotopy equivalences of towers 
$$ \mu: P^{K\vee L}(X) \xra{\sim} P^K(X)\sm P^L(X)$$
with the following properties.
\begin{enumerate}
\item There is a commutative diagram of weak natural transformations:
$$
\xymatrix{
\Sinfty \MapT(K \vee L,X)_+    \ar@{=}[d] \ar[rr]^-{e^{K \vee L}(X)} & &
P^{K \vee L}(X)   \ar[d]_{\mu}^{\wr}  \\
\Sinfty \MapT(K,X)_+ \sm \MapT(L,X)_+ \ar[rr]^-{e^K(X) \sm e^{L}(X)} & &
P^K(X) \sm P^{L}( X). 
}
$$

\item  The induced weak equivalence on $k^{th}$ fibers
$$ F_k^{K\vee L}(X) \xra{\sim} \prod_{i+j = k} F_i^K(X) \sm F_j^L(X)$$
\end{enumerate}
is naturally equivalent  to the product, over $i+j = k$, of the weak natural transformations 
\begin{equation*}
\begin{split}
\Map_{\Sp}((K \vee L)^{( k)}, X^{\sm k})_{h\Sigma_k} & 
\xra{Tr} \Map_{\Sp}((K \vee L)^{( k)}, X^{\sm k})_{h(\Sigma_i \times \Sigma_j)} \\
& \xra{\iota^*} \Map_{\Sp}(K^{(i)}\sm L^{(j)}, X^{\sm k})_{h(\Sigma_i \times \Sigma_j)} \\
& \xla{\sim} \Map_{\Sp}(K^{(i)},  X^{\sm i})_{h\Sigma_i}\sm \Map_{\Sp}(L^{(j)},  X^{\sm j})_{h\Sigma_j},
\end{split} 
\end{equation*}
where $Tr$ is the transfer associated to $\Sigma_i \times \Sigma_j \subset \Sigma_k$, and 
$$\iota:K^{(i)}\sm L^{(j)} \hra (K \vee L)^{(k)}$$
 is the $\Sigma_i \times \Sigma_j$--equivariant inclusion. \\
\end{thm}

Let $\nabla: K \vee K \ra K$ be the fold map.  Since  the diagonal map $\Delta$ is the composite
\begin{equation*}
 \MapT(K,X) \xra{\nabla^*} \MapT(K \vee K,X) = \MapT(K,X) \times \MapT(K,X),  
\end{equation*}
our Product Theorem has consequences for $\Delta$.

Let $\Psi: P^K(X) \ra P^K(X) \sm P^K(X)$ be the weak natural transformation 
$$ P^K(X) \xra{\nabla^*} P^{K \vee K}(X) \xra{\mu} P^K(X) \sm P^K(X).$$

\begin{thm} \label{diag thm} The weak natural transformation $\Psi$ has the following properties.
\begin{enumerate}
\item There is a commutative diagram of weak natural transformations:
$$
\xymatrix{
\Sinfty \MapT(K,X)_+    \ar[d]^{\Sinfty \Delta_+} \ar[rr]^-{e^{K}(X)}  & &
P^{K}(X)   \ar[d]^{\Psi}  \\
\Sinfty \MapT(K,X)_+ \sm \MapT(K,X)_+ \ar[rr]^-{e^K(X) \sm e^{K}(X)} & &
P^K(X) \sm P^{K}( X). 
}
$$
\item The induced weak natural transformation on $k^{th}$ fibers
$$ F^K_k(X) \ra \prod_{i+j=k} F^K_i(X) \sm F_j^K(X)$$
\end{enumerate}
is naturally equivalent to the product, over $i+j = k$, of the composites of the weak natural transformations
\begin{equation*}
\begin{split}
\MapS(K^{(k)}, X^{\sm k})_{h\Sigma_k} & \xra{Tr}
\MapS(K^{(k)}, X^{\sm k})_{h(\Sigma_i \times \Sigma_j)} \\
 & \xra{\pi^*} \MapS(K^{(i)}\sm K^{(j)}, X^{\sm k})_{h(\Sigma_i \times \Sigma_j)} \\ 
& \xla{\sim} \MapS(K^{(i)}, X^{\sm i})_{h\Sigma_i}\sm \MapS(K^{(j)}, X^{\sm j})_{h\Sigma_j}, 
\end{split}
\end{equation*}
where $\pi:K^{(i)} \sm K^{(j)} \ra K^{(k)}$ is the projection. \\
\end{thm}

A typical computational consequence of this would be the following. \\

\begin{cor}  Let $h^*$ be a generalized cohomology theory with products.  Then the associated spectral sequence for computing $h^*(\MapT(K, X))$ is a spectral sequence of bigraded differential graded $h^*$--algebras.  The product on $E_{\infty}^{*,*}$ corresponds to the cup product in $h^*(\MapT(K, X))$ in the usual way, and the product $E_1^{-i,*}\otimes E_1^{-j,*} \ra E_1^{-(i+j),*}$ is induced by the maps on fibers as given in the theorem. \\
\end{cor}

Specializing to $K = S^n$, we have some simplification. \\

\begin{cor}  \label{Sn diagonal cor} There is a natural map of towers
$$ \Psi: P^{S^n}(X) \ra P^{S^n}(X) \sm P^{S^n}(X) $$
under 
$$\Sinfty \Delta_+: \Sinfty (\Omega^n X)_+ \ra 
\Sinfty (\Omega^{n}X \times \Omega^{n}X)_+,$$ such that the associated map on $k^{th}$ fibers is equivalent to the product, over $i+j = k$, of the composites 
\begin{equation*}
\begin{split}
 \C(n,k)_+ \sm_{\Sigma_{k}} (\Sigma^{-n}X)^{\sm k}& \xra{Tr}
\C(n,k)_+ \sm_{\Sigma_i\times \Sigma_j} (\Sigma^{-n}X)^{\sm k} \\
& \ra (\C(n,i) \times \C(n,j))_+ \sm_{\Sigma_i\times \Sigma_j} (\Sigma^{-n}X)^{\sm k} , 
\end{split}
\end{equation*}
where the second map is induced by the $\Sigma_i\times \Sigma_j$--equivariant inclusion 
$$ \C(n,k) \hra \C(n,i) \times \C(n,j).$$
\end{cor}

In this corollary, the second map is an equivalence if $n=\infty$. When $n=1$,  
the $(i,j)^{th}$ component of the map on $k^{th}$ fibers is easily seen to be homotopic to the `shuffle coproduct'
$$ X^{\sm k} \ra X^{\sm i} \sm X^{\sm j},$$
the sum of the $k!/i!j!$ permutations that preserve the order of the first $i$ and last $j$ terms.  Note that this induces the usual product on $E_1$ in the classic Eilenberg--Moore spectral sequence.

\subsection{The Evaluation Theorem}

Our next theorem is a consequence of our study of the map of towers associated to the evaluation maps
$$ \epsilon: K \sm \MapT(K \sm L, X) \ra \MapT(L,X).$$
It is convenient to use reduced towers.  Let $\Tilde{P}^K(X)$ be the fiber of the projection $P^K(X) \ra P^K(*)$.  Then, for all $k$, $P^K_k(X)$ is isomorphic to the product of $\Tilde{P}^K_k(X)$ with the sphere spectrum $S$, and $e^K(X)$ induces a natural transformation 
$$ \Tilde{e}^K(X): \Sinfty \MapT(K,X) \ra \Tilde{P}^K(X).$$

\begin{thm}   There are natural maps of towers 
$$\epsilon:K \sm \Tilde{P}^{K \sm L}(X) \ra \Tilde{P}^L(X)$$
 with the following properties. 
\begin{enumerate}
\item There is a commutative diagram of spectra:
$$
\xymatrix{
\Sinfty K \sm \MapT(K \sm L,X)    \ar[d]^{\Sinfty \epsilon} \ar[rr]^-{1_K \sm \Tilde{e}^{K \sm L}(X)} & &
K \sm \Tilde{P}^{K \sm L}(X)   \ar[d]^{\epsilon}  \\
\Sinfty\MapT(L,X) \ar[rr]^-{\Tilde{e}^{L}(X)} & &
\Tilde{P}^{L}( X). 
}
$$
\item The induced map on $k^{th}$ fibers  is naturally equivalent to the   
composite  
\begin{equation*}
\begin{split} K \sm \Map_{\Sp}((K \sm L)^{( k)}, X^{\sm k})_{h\Sigma_k} & \xra{d^*} K \sm \Map_{\Sp}(K \sm L^{( k)}, X^{\sm k})_{h\Sigma_k} \\
& \xra{\epsilon} \Map_{\Sp}(L^{( k)}, X^{\sm k})_{h\Sigma_k},
\end{split}
\end{equation*}
\end{enumerate}
where the first map is induced by the $\Sigma_k$--equivariant map of spaces
$$ d: K \sm L^{(k)} \ra (K \sm L)^{(k)}$$
which arises by embedding $K$ diagonally in $K^{\sm k}$. \\
\end{thm}

In the 1982 paper \cite{k1}, which studied how the Snaith stable decomposition of $\Omega^n \Sigma^n Y$ interacted with evaluation maps, the second author made use of certain Thom--Pontryagin collapse maps essentially introduced in \cite{may}. These are explicit $\Sigma_k$--equivariant maps of spaces
$$ \beta(m,n,k): S^m \sm \C(m+n,k)_+  \ra S^{mk} \sm \C(n,k)_+.$$

\begin{cor}  \label{Sn evaluation cor} 
There is a natural map of towers
$ \epsilon: \Sigma^m  \Tilde{P}^{S^{m+n}}(X) \ra  \Tilde{P}^{S^n}(X) $ under the evaluation  $\Sinfty \epsilon_+: \Sinfty (\Sigma^m \Omega^{m+n}X)_+ \ra 
\Sinfty (\Omega^n X)_+$, such that the associated map on $k^{th}$ fibers is equivalent to the map
$$ S^m \sm \C(m+n,k)_+ \sm_{\Sigma_k} (\Sigma^{-m-n}X)^{\sm k} \ra \C(n,k)_+ \sm_{\Sigma_k} (\Sigma^{-n}X)^{\sm k}$$
induced by $\beta(m,n,k)$. \\
\end{cor}

We note that the effect in mod $p$ homology of this map on fibers is known, so this theorem can be used computationally. 

\subsection{The $\C(n)$ operad stucture on $P^{S^n}(X)$.}

Our final theorem shows that the little $n$--cubes operad action on $\Omega^n X$ induces an action on our towers in the expected way.

Recall \cite{may} that this action is given by suitably compatible maps
$$ \theta(r): \C(n,r) \times_{\Sigma_r} (\Omega^n X)^r \ra \Omega^nX.$$
Note that $(\Omega^n X)^r = \MapT(\bigvee_r S^n, X)$.  We have the following theorem, which will be made more precise in \secref{operad section}.

\begin{thm}  \label{operad theorem}
For all $n$ and $r$, there is a natural map of towers
$$ \theta(r):\C(n,r)_+ \sm_{\Sigma_r} P^{\bigvee_r S^n}(X) \ra P^{S^n}(X)$$ 
with the following properties.
\begin{enumerate}
\item There is a commutative diagram of spectra:
$$
\xymatrix{
\Sinfty ((\C(n,r) \times_{\Sigma_r} (\Omega^n X)^r)_+)    \ar[d]^{\Sinfty \theta(r)_+} \ar[rr]^-{1 \sm e^{\bigvee_r S^n}(X)} & &
\C(n,r)_+ \sm_{\Sigma_r} P^{\bigvee_r S^n}(X)   \ar[d]^{\theta(r)}  \\
\Sinfty(\Omega^n X)_+ \ar[rr]^-{e^{S^n}(X)} & &
P^{S^n}( X). 
}
$$
\item The associated map on $k^{th}$ fibers is induced by the operad structure
maps
$$ \C(n,r) \times \C(n,k_1) \times \dots \times \C(n,k_r)
\ra \C(n,k),$$
with $k_1 + \dots + k_r = k$.  \\
\end{enumerate}
\end{thm}

Computationally, this implies that the associated spectral sequences for computing mod p homology admit Dyer--Lashof operations.\footnote{Exactly what this statement means is still a matter of investigation by the authors.} \\

\subsection{Organization of the paper.}

The organization of the paper is as follows.  In section 2, we discuss the categories of spectra we work in, and various `naive' constructions including versions of transfer and norm maps.  In section 3, we recall the construction of the Arone tower for $\Sigma^{\infty}\MapT(K,X)_+$, and its homotopical analysis.   We use this in section 4 to prove our Smashing and Evaluation Theorems.  The Product and Diagonal Theorems are proved section 5, after a brief analysis of the smash product of towers.  In section 6 we describe the compatibility  among the various transformations of towers defined in our main theorems.  In section 7, we deduce our various corollaries for the towers $P^{S^n}$, using our explicit equivariant S--duality maps.  Using related constructions with little cubes, \thmref{operad theorem} is proved in section 8, and, in an appendix, a simplified proof of Arone's convergence theorem is given in the case when $K = S^n$.

This paper includes results from the first author's Ph.D. thesis \cite{ahearn}.  The authors wish to thank Greg Arone, Bill Dwyer, and Gaunce Lewis for enlightening mathematical discussions on aspects of this project.  

This research was partially supported by the National Science Foundation.

\section{Background material on spectra} \label{background}

Here we define and discuss various general constructions with spectra that we will later need.  By introducing a small amount of fussiness concerning different universes, all constructions are of a `naive' nature. The material is essentially background, and certainly variations of everything we prove here are already known. 

\subsection{Spectra and universes}

Firstly, we need to specify what we mean by spectra. We find it easiest to work with coordinate free spectra (as in the first pages of \cite{lmms}).  We briefly review the definitions that we need.

Let $\T$ denote the category of compactly generated based spaces. Fixing an infinite dimensional real inner product space $\U$, one defines an associated category of spectra $\Sp\U$.  

An object $X \in \Sp\U$ assigns a space $X(V)$ to every finite dimensional subspace $V \subset \U$, and assigns a structure map $X(V) \ra \Omega^{W-V}X(W)$ to every inclusion $V \subset W$.  Here $W-V$ is the orthogonal complement of $V$ in $W$, and $\Omega^U K = \MapT(S^U,K)$ where $S^U$ is the one point compactification of $U$.  The structure maps are required to be homeomorphisms. 

A map of spectra $f:X \ra Y$ is a collection of maps $f(V): X(V) \ra Y(V)$ compatible with the structure maps in the usual way. This makes $\Sp\U$ into a topological category.

If one deletes the requirement that the structure maps be homeomorphisms, one obtains the category of prespectra $\PP\Sp\U$, and there is a `spectrification' functor $l: \PP\Sp\U \ra \Sp\U$, left adjoint to the inclusion of $\Sp\U$ in $\PP\Sp\U$.  The category $\Sp\U$ has limits and colimits, with limits being formed in $\PP\Sp\U$, and colimits being formed by applying $l$ to the colimit in $\PP\Sp\U$.   When the universe $\U$ is understood and the meaning is clear, we will abbreviate $\Sp\U$ to $\Sp$.

With the elementary constructions to be reviewed later in this section, one does homotopy in the usual way.  The stable category $h\Sp \U$ is then the category obtained from $\Sp \U$ by inverting the weak homotopy equivalences.  A key observation in this approach to spectra is that {\em any} linear isometry $\U \ra \U^{\prime}$ induces the {\em same} equivalence $h\Sp \U  \simeq h\Sp \U^{\prime}$ on passage to homotopy.  We note also that these canonical equivalences are compatible with the various constructions given below. See \cite[chapter II]{lmms} for more detail.

\subsection{Suspension spectra}

There is an adjoint pair
$$ \T 
\begin{array}{c} \stackrel{\Sinfty}{\longleftarrow} \\ \stackrel{\longrightarrow}{\Omega^{\infty}} 
\end{array}  
\Sp\U   $$
defined by $\Omega^{\infty}X = X(0)$ and $(\Sinfty K)(V) = \colim_W \Omega^W \Sigma^{W + V}K$. Here $\Sigma^U K$ denotes $S^U \sm K$, as usual. We let $Q = \Omega^{\infty}\Sinfty: \T \ra \T$, and $S = \Sinfty S^0$.  

When it is necessary to remember $\U$, we will use the notation $\Sinfty_{\U}$, etc. (This follows our general rule with all constructions involving spectra: we will be notationally pedantic when it seems prudent.)

\subsection{Stablization, elementary smash products and function spectra} \label{smash products of spectra}

Given $K \in \T$ and $X \in \Sp\U$, we define spectra $K \sm X$ and $\Map_{\Sp\U}(K,X)$ in $\Sp\U$ as follows: $K \sm X$ is the spectrification of the prespectrum with $V^{th}$ space $K \sm X(V)$, and $\Map_{\Sp\U}(K,X)(V) = \MapT(K,X(V))$.  These constructions are adjoint to each other,
$$ \Hom_{\Sp\U}(K\sm X,Y) = \Hom_{\Sp\U}(X, \Map_{\Sp\U}(K,Y)),$$
and one can deduce various useful isomorphisms in $\Sp\U$ \  \cite[p.17, p.20]{lmms}:
$$ \Map_{\Sp\U}(K\sm L,X) = \Map_{\Sp\U}(K,\Map_{\Sp\U}(L,X)),$$
$$ (K\sm L) \sm X = K \sm (L\sm X),$$
and 
$$ K \sm \Sinfty L = \Sinfty (K \sm L).$$
When clear from context, we will write $\Map_{\Sp \U}(K,X)$ as $\MapS(K,X)$, and then $\MapS(K,\Sinfty L)$ as $\MapS(K,L)$.

A `stabilization' map
$$ s: \Sinfty \MapT(K,L) \ra \MapS(K, L)$$
can now be defined as the adjoint to 
$$ \MapT(K,L) \xra{\MapT(K,\eta)} \MapT(K, QL) = \Omega^{\infty}\MapS(K,L)$$
where $\eta: L \ra QL$ is adjoint to the identity on $\Sinfty L$.  Also arising from adjunctions are evaluation maps
$$\epsilon_{\T}: K \sm \MapT(K,L)  \ra L$$
and 
$$\epsilon_{\Sp}: K \sm \MapS(K,X) \ra X.$$
The next lemma is proved with formal categorical arguments.

\begin{lem} \label{eval lemma}  For any spaces $K$ and $X$, there is a commutative diagram
$$
\xymatrix{
K \sm \Sinfty \MapT(K,X) \ar@{=}[d] \ar[r]^-{1 \sm s} &
K \sm \MapS(K, X) \ar[d]^{\epsilon_{\Sp}}  \\
\Sinfty(K \sm \MapT(K,X)) \ar[r]^-{\Sinfty \epsilon_{\T}} &
\Sinfty X. 
}
$$
\end{lem}

A variation on these constructions goes as follows. (Compare with \cite[pp.68,69]{lmms}.)
Given two universes $\U$ and $\U^{\prime}$, there is an external smash product
$$ \sm: \Sp\U \times \Sp\U^{\prime} \ra \Sp(\U\oplus\U^{\prime})$$
defined by letting $X\sm Y$ be the spectrification of the prespectrum with $(V\oplus W)^{th}$ space $X(V) \sm Y(W)$.\footnote{This formula suffices because subspaces of $\U \oplus \U^{\prime}$ of the form $V \oplus W$ are cofinal among all finite dimensional subspaces.}  Dually, there is an external mapping spectrum functor:
$$ \Map_{\Sp\U}: \Sp\U^{op} \times \Sp(\U \oplus \U^{\prime}) \ra \Sp\U^{\prime}$$
defined by $\Map_{\Sp\U}(X, Z)(W) = \Hom_{\Sp\U}(X,Z_W)$, where $Z_W(V) = Z(V\oplus W)$.  

Again these constructions are adjoint:
$$ \Hom_{\Sp(\U\oplus \U^{\prime})}(X \sm Y, Z) = \Hom_{\Sp \U}(X, \Map_{\Sp \U^{\prime}}(Y,Z)).$$
Again one can formally deduce useful properties, e.g.\ there are isomorphisms in $\Sp(\U\oplus \U^{\prime})$:
$$  \Sinfty_{\U}K \sm \Sinfty_{\U^{\prime}}L = \Sinfty_{\U \oplus \U^{\prime}} (K \sm L).$$
Another useful property, which follows by a check of the definitions, is that, for all $K \in \T$ and $Z \in \Sp(\U\oplus \U^{\prime})$, there are isomorphisms in $\Sp(\U^{\prime})$:
$$ \Map_{\Sp \U}(\Sinfty_{\U}K, Z) = \Map_{\Sp \U^{\prime}}(K, i^*Z),$$
where $i: \U^{\prime} \ra  \U\oplus \U^{\prime}$ is the inclusion.

Given $X \in \Sp \U$ and $Z\in \Sp(\U\oplus \U^{\prime})$, there is an evaluation map
$$ \epsilon: X \sm \Map_{\Sp \U}(X, Z) \ra Z.$$
Precomposing this with $\epsilon: K \sm \Map_{\Sp \U}(K,X) \ra X$ and then adjointing defines a composition map
$$  \circ: \Map_{\Sp \U}(K,X) \sm \Map_{\Sp \U}(X,Z)\ra \Map_{\Sp(\U\oplus \U^{\prime})}(K,Z)$$
for all $K \in \T$.  We will use this construction when defining norm maps in \secref{norm maps} below.

Given spaces $K$, $L$, and spectra $X \in \Sp\U$, $Y \in \Sp\U^{\prime}$, we define
$$ \sm: \Map_{\Sp\U}(K, X) \sm \Map_{\Sp\U^{\prime}}( L, Y) \ra \Map_{\Sp(\U \oplus \U^{\prime})}(K \sm L, X \sm Y)$$
to be adjoint to the composite of the natural isomorphism
{\small
$$ (K \sm L) \sm \Map_{\Sp\U}(K, X) \sm \Map_{\Sp\U^{\prime}}( L, Y) = 
(K  \sm \Map_{\Sp\U}(K, X)) \sm (L \sm \Map_{\Sp\U^{\prime}}( L, Y)) $$
with \  ${e_{\Sp}\sm e_{\Sp}}:(K  \sm \Map_{\Sp\U}(K, X)) \sm (L \sm \Map_{\Sp\U^{\prime}}( L, Y)) \ra X \sm Y.$}

This is analogous to the usual pairing between mapping spaces
$$ \sm: \MapT(K, X) \sm \MapT( L, Y) \ra \MapT(K \sm L, X \sm Y),$$
and the next lemma records that these constructions are compatible under stabilization. 

\begin{lem} \label{smash lemma}  For any spaces $K$, $L$, $X$, $Y$, and universes $\U$, $\U^{\prime}$, there is a commutative diagram in $\Sp(\U \oplus \U^{\prime})$ 
{\small $$
\xymatrix{
\Sinfty_{\U} \MapT(K,X)\sm \Sinfty_{\U^{\prime}} \MapT(L,Y) \ar@{=}[r] \ar[d]^{s \sm s} &
\Sinfty_{\U \oplus \U^{\prime}} (\MapT(K,X)\sm \MapT(L,Y))  \ar[d]^{\Sinfty \sm}    \\
 \Map_{\Sp \U}(K, \Sinfty_{\U}X) \sm \Map_{\Sp \U^{\prime}}( L,  \Sinfty_{\U^{\prime}}Y) \ar[d]^{\sm} &
\Sinfty_{\U \oplus \U^{\prime}} \MapT(K\sm L, X\sm Y)\ar[d]^{s} \\
\Map_{\Sp(\U \oplus \U^{\prime})}(K \sm L, \Sinfty_{\U}X \sm \Sinfty_{\U^{\prime}}Y)  \ar@{=}[r]&    \Map_{\Sp(\U \oplus \U^{\prime})}(K \sm L, \Sinfty_{\U \oplus \U^{\prime}} (X \sm Y))
}
$$}
\end{lem}

Once again, this is proved with formal categorical arguments.

The next lemma is standard.

\begin{lem} \label{smash lemma 2} If $K$ and $L$ are finite CW complexes, then
$$ \sm: \Map_{\Sp\U}(K, X) \sm \Map_{\Sp\U^{\prime}}( L, Y) \ra \Map_{\Sp(\U \oplus \U^{\prime})}(K \sm L, X \sm Y)$$
is a weak homotopy equivalence.
\end{lem}

\subsection{Spaces and spectra of natural transformations}

If $\J$ is a small category, and $K:\J \ra \T$ and $X:\J \ra \Sp$ are two functors of the same variance, we will write $\Map_{\Sp}^{\J}(K,X)$ for the spectrum constructed as the categorical equalizer in $\Sp$ of the two evident maps
$$ \prod_{j \in Ob(\J)} \MapS(K(j), X(j)) 
\begin{array}{c} {\longrightarrow} \\ {\longrightarrow}
\end{array}
\prod_{\alpha:j^{\prime} \ra j^{\prime \prime} \in Mor(\J)} \MapS(K(j^{\prime}), X(j^{\prime \prime})).
$$
Similarly, if $K,X: \J \ra \T$ are two functors of the same variance, one gets a space $\Map_{\T}^{\J}(K,X)$, which can be interpreted as the space of natural transformations from $K$ to $X$. The stable and unstable constructions are related by $\Map_{\Sp}^{\J}(K,X)(V) = \Map_{\T}^{\J}(K,X(V))$, for any $V \in \U$.

It is useful to observe that if $\J = G$, a finite group viewed as a category with one object, then  $\Map_{\Sp}^{G}(K,X)$ is the categorical fixed point spectrum of the naive $G$--spectrum $\MapS(K,X)$ with conjugation $G$--action.  In this case, we also write $X/G$ for the categorical orbit spectrum.

\subsection{Norm maps, transfers, and Adams isomorphisms} \label{norm maps}

In this subsection, we give quick definitions of transfer and norm maps suitable for our later homotopical identification of natural transformations between fibers in the Arone towers.   These definitions are adapted to our setting, but are intended to agree in the homotopy category with anyone else's transfer and norm maps.  As far as the authors can tell, constructions of norm maps using only ``naive'' constructions first appeared in the literature in the 1989 paper of Weiss and Williams \cite[\S 2]{ww}. (Those authors credit Dwyer with some of these ideas, and, of course, Adams' paper \cite{adams} was influential.)  Our definitions are small perturbations of those in the recent preprint of John Klein \cite{klein}.  \propref{transfer prop}, which relates transfer and norm maps, appears to be new in the literature, and a desire for a transparent proof of this has guided our constructions.

Let $G$ be a finite group, and call a spectrum with $G$--action a $G$--spectrum.  Fix two universes $\U$ and $\U^{\prime}$, and let $i: \U^{\prime} \ra  \U\oplus \U^{\prime}$ be the inclusion.

\begin{defns}  Given a subgroup $H \leq G$, and a $G$--spectrum $X \in \Sp(\U\oplus \U^{\prime})$, we define the homotopy fixed point and homotopy orbit spectra as follows.
\begin{enumerate}
\item  $X^{hH} = \Map_{\Sp(\U\oplus \U^{\prime})}^H(EG_+, X).$
\item  $X_{hH} = (EG_+ \sm \Map_{\Sp \U}^H(EG_+, \Sinfty_{\U} G_+)\sm i^*X)/G$.
\end{enumerate}
\end{defns}

The first definition is, we trust, expected.  The corollary of the next lemma says that the second has the correct homotopy type.

\begin{lem} \label{G/H lem} There is a weak equivalence of $G$--spectra in $\Sp \U$
$$ \Sinfty G/H_+ \simeq \Map^H_{\Sp}(EG_+, G_+).$$
\end{lem}
\begin{proof}  There are weak equivalences and isomorphisms of $G$--spectra:
\begin{equation*}
\begin{split}
\Sinfty G/H_+ & \xra{\sim} \MapS(G/H_+, S) \\
  & =\Map_{\Sp}^H(G_+, S) \\
 & \xra{\sim} \Map_{\Sp}^H(EG_+ \sm G_+, S) \\
 & = \Map_{\Sp}^H(EG_+, \Map_{\Sp}(G_+, S)) \\
 & \xla{\sim} \Map_{\Sp}^H(EG_+, G_+).
\end{split}
\end{equation*}
Here the first and last maps arise in the same manner.  If $K < G$ is any subgroup, there is a commutative diagram of $G$--spectra
\begin{equation*}
\xymatrix{
 \bigvee_{gK \in G/K}S \ar[d]^{\simeq} \ar@{=}[r] &
\Sinfty G/K_+ \ar[d]  \\
 \prod_{gK \in G/K}S \ar@{=}[r] &
\MapS(G/K_+, S)
}
\end{equation*}
where the left vertical map is the inclusion of the wedge into the product, a weak homotopy equivalence.
\end{proof}

\begin{cor} \label{hoorbit cor} There is a weak natural equivalence in $\Sp(\U\oplus \U^{\prime})$
$$X_{hH} \simeq (EG_+ \sm X)/H.$$
\end{cor}
\begin{proof}  There are weak natural equivalences  
$$X_{hH} \simeq (EG_+ \sm \Sinfty_{\U} G/H_+ \sm i^* X)/G \xra{\sim} (EG_+ \sm  G/H_+ \sm X)/G = (EG_+ \sm X)/H.$$
Here the first equivalence is a consequence of the lemma, and the second follows from the fact that, very generally, there is a natural weak equivalence $\Sinfty_{\U} K \sm i^* X \ra K \sm X$.
\end{proof}

Our transfer maps are defined as follows.

\begin{defns}  Let $K \leq H \leq G$, and let $X$ be a $G$--spectrum in $\Sp(\U \oplus \U^{\prime})$.
\begin{enumerate}
\item Let $tr_K^H: \Map^H_{\Sp \U}(EG_+, G_+) \ra \Map^K_{\Sp \U}(EG_+, G_+)$ be the inclusion of fixed point spectra.
\item Let $Tr_K^H(X): X_{hH} \ra X_{hK}$ be the natural map induced by $tr_K^H$. \\
\end{enumerate}
\end{defns}

We sketch a proof that $Tr_K^H(X)$, viewed as a natural transformation of functors on the homotopy category of spectra with $G$--action, agrees with other standard constructions of the transfer, in particular, the transfer arising from \cite{lmms}.  Both of these transfers behave well with respect to pushouts and weak equivalences in the $X$ variable, and with respect to forgetful functors arising from subgroup inclusions.  Using these facts one can reduce to just needing to show that the two definitions of $Tr_H^G(G_+)$ agree up to weak equivariant homotopy.  For us, this map is equivalent to the map
$$ S \xra{\sim} \MapS(G/G_+, S) \ra \MapS(G/H_+,S)$$
induced by the projection $\pi: G/H_+ \ra G/G_+$.  Now one checks that this agrees with the composite
$$ S \xra{\tau} \Sinfty G/H_+ \xra{\sim} \MapS(G/H_+, S)$$
where $\tau$ is the pretransfer of \cite[p.181]{lmms}.  \\

We now define our norm maps.

\begin{defn}  Given $H \leq G$, and a $G$--spectrum $X$ in $\Sp(\U \oplus \U^{\prime})$, let 
$$\Phi_H(X): X_{hH} \ra X^{hH}$$
 be defined as follows.  First note that  
$$i^*X = \Map_{\Sp \U^{\prime}}^G(G_+, i^* X) = \Map_{\Sp \U}^G(\Sinfty_{ \U}G_+, X) .$$
Now consider composition
$$  \circ: \Map_{\Sp \U}(EG_+,\Sinfty_{ \U}G_+) \sm \Map_{\Sp \U}^G(\Sinfty_{\U}G_+,X)\ra \Map_{\Sp(\U\oplus \U^{\prime})}(EG_+,X).$$
This is $G$--equivariant with respect to the usual conjugation $G$--action on the two terms $\Map_{\Sp \U}(EG_+,\Sinfty_{ \U}G_+)$ and $\Map_{\Sp(\U\oplus \U^{\prime})}(EG_+,X)$.  Taking $H$ fixed points then yields a map of spectra
$$ \Map_{\Sp \U}^H(EG_+,\Sinfty_{\Sp \U}G_+) \sm i^*X  \ra X^{hH}.$$ 
Now one notes that this map is invariant with respect to the diagonal $G$--action on the domain, where $G$ acts on $(\Map_{\Sp \U}^H(EG_+,\Sinfty_{\Sp \U}G_+)$ by acting on the right on $G_+$.  Thus one has an induced map
$$ (\Map_{\Sp \U}^H(EG_+,\Sinfty_{\Sp \U}G_+) \sm i^*X)/G  \ra X^{hH}.$$ 
$\Phi_H(X)$ is then obtained by precomposing this map with the map 
$$ X_{hH} \ra (\Map_{\Sp \U}^H(EG_+,\Sinfty_{\Sp \U}G_+) \sm i^*X)/G $$
induced by $EG_+ \ra S^0$. \\
\end{defn}

By construction, the following proposition is self evident.

\begin{prop} \label{transfer prop} Given $K \leq H \leq G$, and $X \in \Sp(\U \oplus \U^{\prime})$, there is a commutative diagram of spectra
\begin{equation*}
\xymatrix{
X_{hH} \ar[d]^{Tr_K^H(X)} \ar[r]^{\Phi_H(X)} &
X^{hH} \ar[d]  \\
X_{hK} \ar[r]^{\Phi_K(X)} &
X^{hK}
}
\end{equation*}
where the unlabelled vertical arrow is the inclusion of fixed point spectra.
\end{prop}

Norm maps should be equivalences under suitable freeness and finiteness conditions.  In homotopy, such equivalences have been termed `Adams isomorphisms' \cite{lmms}.  The version we need goes as follows.

\begin{prop} \label{Adams iso}  If $X = \MapS(K,Y)$, where Y is any $G$--spectrum and $K$ is any finite free $G$--CW complex, then the norm map
$$ \Phi_H(X): X_{hH} \ra X^{hH}$$
is a weak homotopy equivalence for all $H<G$.
\end{prop}
\begin{proof}  We show that $\Phi_H(X)$ is an equivalence in various cases.  When $X = \Sinfty_{\U \oplus \U^{\prime}}G_+$, via the weak equivalences of \lemref{G/H lem}, $\Phi_H(X)$ corresponds to the isomorphism
$$ (\Sinfty_{\U}G/H_+ \sm \Sinfty_{\U^{\prime}}G_+)/G = \Sinfty_{\U \oplus \U^{\prime}} G/H_+.$$
Now we note that both the domain and range of $\Phi_H(X)$ preserve equivariant weak homotopy equivalences and homotopy cofiber sequences.  Thus, by induction on cells, $\Phi_H(X)$ is an equivalence if $Y$ is any $G$--spectrum equivalent to a finite $G$--CW spectrum, and $K$ is any finite free $G$--CW complex.  Now we note that, if $K$ is a finite free $G$--CW complex, then both the domain and range of $\Phi_H(\MapS(K,Y))$ commute with homotopy colimits in the $Y$ variable.  The proposition follows.
\end{proof}

An application of this that we will need later goes as follows.

\begin{cor} \label{smash cor} Let $G$ and $H$ be two finite groups.  If $X$ is a $G$--spectrum, $K$ a finite free $G$--CW complex, $Y$ an $H$--spectrum, and $L$ a finite free $H$--CW spectrum, then 
$$ \sm: \Map_{\Sp\U}^G(K, X) \sm \Map_{\Sp\U^{\prime}}^H( L, Y) \ra \Map_{\Sp(\U \oplus \U^{\prime})}^{G \times H}(K \sm L, X \sm Y)$$
is a weak homotopy equivalence.
\end{cor}
\begin{proof}  Let $A=\Map_{\Sp}(K, X)$, $B=\Map_{\Sp}( L, Y)$, and $C=\Map_{\Sp}(K \sm L, X \sm Y)$.  We wish to show that the $G \times H$--equivariant weak equivalence $A \sm B \ra C$ induces an equivalence $A^G \sm B^H \ra C^{G \times H}$.  Since there are equivalences $A^G \xra{\sim} A^{hG} \xla{\sim} A_{hG}$, and similarly for $B$ and $C$, it suffices to show that $A_{hG} \sm B_{hH} \ra C_{h(G\times H)}$ is an equivalence.  But this is clear, as easy formal arguments show that $A_{hG} \sm B_{hH} = (A \sm B)_{h(G \times H)}$.
\end{proof}

Finally we note:

\begin{cor} \label{norm cor} Let $G$ be a finite group.  If $K$ is a finite free $G$--CW complex, then for all $G$--spectra $Y \in \Sp (\U \oplus \U^{\prime})$, there is a weak natural equivalence
$$ (\Map_{\Sp \U}(K,S) \sm i^*Y)_{hG} \simeq \Map_{\Sp (\U \oplus \U^{\prime})}^G(K,Y)$$
\end{cor}
\proof  This arises from the equivalences
$$ (\MapS(K,S) \sm Y)_{hG} \xra{\sim} \MapS(K,Y)^{hG} \xla{\sim} \Map_{\Sp}^G(K,Y)\eqno{\qed}.$$

\section{The Arone model} \label{arone model}
In this section we review the Greg Arone's explicit construction of the tower $P^K(X)$ of fibrations of spectra under $\Sinfty \MapT(K,X)_+$, tweaked a bit to lend itself to the homotopical analysis we are interested in.\footnote{For example, Arone works with the functors from (based) spaces to spaces which send a space $X$ to $Q\MapT(K,X)$.  It seems more natural to regard the basic functors as going from spaces to a suitable category of spectra.} 

\subsection{Definition of the tower}

We introduce a small category central to our work. \\

\begin{defn}  Let $\E$ denote the category with objects the finite sets $\mathbf  0 = \emptyset$ and  ${\mathbf  n} = \{1,2,\dots, n\}$, $n \geq 1$, and morphisms the surjective functions. This has full subcategories $\Tilde \E$ whose objects are $\mathbf n$ with $n \geq 1$, $\E_k$, whose objects are $\mathbf n$ with $n \leq k$, and $\Tilde \E_k = \Tilde \E \cap \E_k$.\footnote{This category $\E$ of epimorphisms was called $\mathcal M$ by Arone in \cite{arone} and $\mathcal M^{op}$ by McCarthy in \cite{mccarthy2}.  The senior author objects.}
\end{defn}

Fundamental $\E$--spaces are the following.

\begin{defn}  If $X$ is a pointed space, let $X^{\sm}: \E^{op} \ra \T$ be the following functor.  On objects, it assigns to ${\mathbf n}$, $X^{\sm n}$, the $n$--fold smash product of $X$ with itself, where we use the convention that $X^{\sm 0} =S^0$.  On morphisms, it assigns to a surjective function $\alpha: {\mathbf  n} \ra {\mathbf  m}$, the associated `diagonal' map $\alpha^*:X^{\sm m}\ra X^{\sm n}$ sending $x_1 \sm \dots \sm x_m$ to $x_{\alpha(1)} \sm \dots \sm x_{\alpha(n)}$. Note that this is well defined precisely because $\alpha$ is surjective. \\
\end{defn}

Armed with these $\E$--spaces, we can define Arone's towers in our setting. \\

\begin{defns}  Given two spaces $K$ and $X$ in $\T$, define spectra $P_{\infty}^K(X)$, $P^K_k(X)$, $\Tilde{P}^K_{\infty}(X)$, and $\Tilde{P}_k^K(X)$, by the formulae 
$$P_{\infty}^K(X) = \Map_{\Sp}^{\E}(K^{\sm}, X^{\sm}),$$   
$$P_k^K(X) = \Map_{\Sp}^{\E_k}(K^{\sm}, X^{\sm}),$$  
$$\Tilde{P}^K_{\infty}(X) = \Map_{\Sp}^{\Tilde \E}(K^{\sm}, X^{\sm}),$$ 
$$\Tilde{P}_k^K(X) = \Map_{\Sp}^{\Tilde \E_k}(K^{\sm}, X^{\sm}).$$
\end{defns}

There are evident restriction maps
$$ p_k: P^K_{k}(X) \ra P^K_{k-1}(X),$$
defining a tower $P^K(X)$, compatible with restrictions maps
$$q_k: P^K_{\infty}(X) \ra P^K_k(X),$$
and $ P^K_{\infty}(X) = \lim_{k} P^K_k(X).$
As we will discuss below, the maps $p_k$ are fibrations in $\Sp$; thus $P^K_{\infty}(X)$ is also equivalent to the homotopy limit of the tower. \\

\begin{notation} Let $F^K_k(X)$ be the fiber of $p_k:P^K_k(X) \ra P^K_{k-1}(X)$. \\
\end{notation}

The reduced functors are related to the unreduced functors by 
$$ P^K_k(X) = \Tilde{P}_k^K(X) \times S,$$
(with the product in $\Sp$), and they similarly form a tower with limit (and holimit) $\Tilde{P}^K_{\infty}(X)$. \\

Now we define natural transformations
$\Sinfty \MapT(K,X)_+ \ra P^K(X).$ \\

For $n \geq 0$, let $$e^K(X,{n}): \Sinfty \MapT(K,X)_+ \ra \Map_S(K^{\sm n}, X^{\sm n})$$ be the composite
$$\Sinfty \MapT(K,X)_+ \xra{\Sinfty \xi_n} \Sinfty \MapT(K^{\sm n}, X^{\sm n}) \xra{s} \MapS(K^{\sm n}, X^{\sm n}),$$
where 
$$ \xi_n: \MapT(K,X)_+ \ra \MapT(K^{\sm n}, X^{\sm n})$$
sends a function $f$ to $f^{\sm n}$. (When $n=0$, this means the evident projection $\MapT(K,X)_+ \ra S^0$.)  \\

\begin{lem}  Let $\alpha:{\mathbf n}\ra {\mathbf m}$ be a surjective function.  Then, for all $K$ and $X$ there is a commutative diagram in $\Sp$
$$
\xymatrix{
\Sinfty \MapT(K,X)_+    \ar[d]^{e^K(X,m)} \ar[r]^{e^K(X,n)} &
\MapS(K^{\sm n}, X^{\sm n})  \ar[d]^{\alpha_*}  \\
\MapS(K^{\sm m}, X^{\sm m}) \ar[r]^{\alpha^*} &
\MapS(K^{\sm m}, X^{\sm n}). 
}
$$
\end{lem}

\begin{proof}  The analogous diagram in $\T$,
$$
\xymatrix{
 \MapT(K,X)_+    \ar[d]^{\xi_m} \ar[r]^-{\xi_n} &
\MapT(K^{\sm n}, X^{\sm n})  \ar[d]^{\alpha_*}  \\
\MapT(K^{\sm m}, X^{\sm m}) \ar[r]^{\alpha^*} &
\MapT(K^{\sm m}, X^{\sm n}), 
}
$$
clearly commutes, and the lemma follows, using the naturality of $s$.
\end{proof}

This lemma says that the maps 
$$\prod_n  e^K(X,{n}): \Sinfty \MapT(K,X)_+ \ra \prod_n \MapS(K^{\sm n}, X^{\sm n})$$
are equalized by the maps defining $P^K_{\infty}(X) = \Map_{\Sp}^{\E}(K^{\sm}, X^{\sm})$.  Thus they define
$$e^K_{\infty}(X): \Sinfty \MapT(K,X)_+ \ra P^K_{\infty}(X),$$
and then 
$$e^K_k(X) = q_k \circ e^K_{\infty}(X): \Sinfty \MapT(K,X)_+ \ra P^K_k(X). $$
The main theorem of \cite{arone} is a convergence result. In our context, it reads as follows.  \\

\begin{thm} \label{arone thm} Let $K$ be a finite CW complex and $X$ a space with connectivity at least as large as the dimension of $K$.  Then $e^K_k(X)$ is $(1+ \text{conn }X - \text{dim } K)(1+k)-1$ connected.  In particular, $e^K_{\infty}(X)$ is a weak homotopy equivalence and the tower is strongly convergent.
\end{thm}

In \appref{convergence appendix}, we will outline how the proof of this theorem goes in the case when $K$ is a sphere.   Needed first in any proof, however, is an analysis of the homotopical behavior of the tower $P^K_*(X)$.  We now proceed to make this analysis, as these results are also needed to prove our main theorems.

\subsection{Homotopical analysis of the tower}
The key to understanding $P^K(X)$ homotopically is to take advantage of two observations.  Firstly, the filtration of $\E$ by the subcategories $\E_k$ induces a natural filtration on contravariant functors from $\E$ to $\T$ or $\Sp$.  Secondly, this filtration on the particular functors $K^{\sm}$ is particularly nice, and essentially exhibits them as cofibrant $\E^{op}$--objects in the functor categories.

We should say immediately that these sorts of observations have been made before.  See, for example, \cite[Appendix A]{mccarthy2}, \cite[\S 3]{lydakis}, as well as Arone's own paper \cite{arone}.  These are all modern references, but, slightly disguised, these ideas are certainly much older.

With $\C$ either $\T$ or $\Sp$, let $\C^{\E}$ denote the category of contravariant functors $X:\E^{op} \ra \C$. The inclusion $i_k: \E_k \hra \E$ induces the restriction $i_k^*: \C^{\E}\ra \C^{\E_k}$ with left adjoint ${i_k}_*: \C^{\E_k}\ra \C^{\E}$, and we let $X_k = {i_k}_*i_k^*X$.  Formally, one sees that, for all $X \in \T^{\E}$ and $Y \in \C^{\E}$, there are natural isomorphisms in $\C$
$$ \Map_{\C}^{\E}(X_k, Y) = \Map_{\C}^{\E_k}(X,Y).$$
Explicitly,
$$ X_k(n) = \colim_{{\mathbf n}\da \E_k} \Tilde{X},$$
where ${\mathbf n}\da \E_k$ denotes the usual category under $\mathbf n$ with objects $\mathbf n \ra \mathbf j$ in $\E$ with $j \leq k$, and $\Tilde{X}(\mathbf n \ra \mathbf j) = X(j)$.

The $X_k$ assemble into a filtration of $X$,
$$ X_0 \ra X_1 \ra X_2 \ra \dots, $$
and, noting that $X_k(n) = X(n)$ for all $k \geq n$, one sees that $X$ is realized as the colimit.

In the next lemma, $\Sigma_k$ denotes the symmetric group on $k$ letters, viewed as the morphisms $\mathbf k \ra \mathbf k$ in $\E$.\\

\begin{lem} For all $X \in \T^{\E}$, there is a pushout in $\T^{\E}$
$$
\xymatrix{
\E( \ , \mathbf k)_+ \sm_{\Sigma_k} X_{k-1}(k)    \ar[d] \ar[r] &
X_{k-1} \ar[d]  \\
\E( \ , \mathbf k)_+ \sm_{\Sigma_k} X_{k}(k) \ar[r] &
X_k. 
}
$$
\end{lem}

\begin{proof}  It is very easy to check that the diagram of small categories
$$
\xymatrix{
\E( \mathbf n , \mathbf k) \times_{\Sigma_k} (\mathbf k \da \E_{k-1})    \ar[d] \ar[r] &
(\mathbf n \da \E_{k-1})  \ar[d]  \\
\E( \mathbf n , \mathbf k) \times_{\Sigma_k} (\mathbf k \da \E_{k})  \ar[r] &
(\mathbf n \da \E_{k})  
}
$$
specializes to a pushout diagram of finite sets on both objects and morphisms.  

It follows that there is a pushout diagram in $\T$
$$
\xymatrix{
\colim_{\E( \mathbf n , \mathbf k) \times_{\Sigma_k} (\mathbf k \da \E_{k-1})} \tilde{X}    \ar[d] \ar[r] &
\colim_{(\mathbf n \da \E_{k-1})} \tilde{X}  \ar[d]  \\
\colim_{\E( \mathbf n , \mathbf k) \times_{\Sigma_k} (\mathbf k \da \E_{k})} \tilde{X}  \ar[r] &
\colim_{(\mathbf n \da \E_{k})} \tilde{X}.  
}
$$
This last diagram rewrites as the diagram of the lemma, evaluated at $\mathbf n$.
\end{proof}

\begin{cor}  For all $X \in \T^{\E}$ and $Y \in \C^{\E}$ there is a pullback diagram in $\C$
$$
\xymatrix{
\Map_{\C}^{\E}(X_k, Y)    \ar[d] \ar[r] & 
\Map_{\C}^{\Sigma_k}(X_k(k), Y(k))  \ar[d]  \\
\Map_{\C}^{\E}(X_{k-1}, Y) \ar[r] &
\Map_{\C}^{\Sigma_k}(X_{k-1}(k), Y(k)). 
}
$$
\end{cor}

Now suppose that $X=K^{\sm}$.  Observe that $K^{\sm}_k(k) = K^{\sm k}$, and that $K^{\sm}_{k-1}(k)$ is the fat diagonal $\Delta_k(K)$ inside $K^{\sm k}$.    If $K$ is a CW complex, then $K^{\sm k}$ can be obtained from $\Delta_k(K)$ by attaching only free $\Sigma_k$--cells.  The fact that $(K^{\sm k},\Delta_k(K))$ is an equivariant CW pair implies that the inclusion $\Delta_k(K) \ra K^{\sm k}$ is an equivariant cofibration.  Recalling that $K^{(k)}$ denotes $K^{\sm k}/\Delta_k(K)$ (as in \cite{arone}), we conclude:

\begin{prop} \label{fiber prop} $\Map_{\Sp}^{\Sigma_k}(K^{\sm k}, X^{\sm k}) \ra \Map_{\Sp}^{\Sigma_k}(\Delta_k(K), X^{\sm k})$ is a fibration with fiber $\Map_{\Sp}^{\Sigma_k}(K^{( k)}, X^{\sm k})$.  Thus $p_k: P^K_k(X) \ra P^K_{k-1}(X)$ is a fibration with fiber $F^K_k(X) = \Map_{\Sp}^{\Sigma_k}(K^{( k)}, X^{\sm k})$.  Furthermore, there are natural weak equivalences
$$ F^K_k(X) \xra{\sim} \Map_{\Sp}(K^{( k)}, X^{\sm k})^{h \Sigma_k} \xla{\sim} (\Map_{\Sp}(K^{( k)}, S) \sm  X^{\sm k})_{h \Sigma_k}.$$
\end{prop}

The last statement here makes it clear that the Arone tower has the form of a Goodwillie tower.   In the language of \cite{goodwillie3}, we conclude:

\begin{cor} The $k^{th}$ Taylor coefficient of the functor sending  a space $X$ to the spectrum $\Sinfty \MapT(K,X)$ is the $\Sigma_k $--spectrum $\MapS(K^{(k)}, S)$.
\end{cor}

\section{The Smashing and Evaluation Theorems} \label{evaluation section}

\subsection{The Smashing Theorem}  The first of our general theorems studying additional structure in the Arone tower involves smashing with a constant space $K$.  

There is an unstable map
$$ \eta: \MapT(L,X) \ra \MapT(K\sm L, K\sm X)$$
which sends a function $f:L \ra X$ to $1_K \sm f: K \sm L \ra K \sm X$.
Note that this map can be written as the composite
$$ \MapT(L,X) \ra \MapT(K\sm L, K\sm L \sm \MapT(L,X)) \ra \MapT(K\sm L, K\sm X)$$
where the first map is a unit of an adjunction, and the second map is induced by $\epsilon_{\T}: L \sm \MapT(L,X) \ra X$.  Replacing $\T$ by $\Sp$ in this composite then similarly defines
$$  \eta: \MapS(L,X) \ra \MapS(K\sm L, K\sm X),$$
and the unstable and stable maps will be compatible under stabilization in the evident way.

\begin{thm} \label{smash theorem}  There are natural maps of towers 
$$\eta:P^L(X) \ra P^{K \sm L}(K \sm X)$$
 with the following properties. 
\begin{enumerate}
\item There is a commutative diagram in $\Sp$:
$$
\xymatrix{
\Sinfty\MapT(L,X)_+    \ar[d]^{\eta} \ar[rr]^-{e^L(X)} & &
P^L(X)   \ar[d]^{\eta}  \\
\Sinfty\MapT(K \sm L,K \sm X)_+ \ar[rr]^-{e^{K \sm L}(K \sm X)} & &
P^{K \sm L}(K\sm X). 
}
$$
\item The induced map on $k^{th}$ fibers is the  composite 
\begin{equation*}
\begin{split}
\Map_{\Sp}^{\Sigma_k}(L^{( k)}, X^{\sm k}) & 
\xra{\eta} \Map_{\Sp}^{\Sigma_k}(K^{\sm k} \sm L^{( k)}, K^{\sm k} \sm X^{\sm k}) \\
& \xra{p^*} \Map_{\Sp}^{\Sigma_k}((K\sm L)^{( k)}, (K \sm X)^{\sm k}),
\end{split} 
\end{equation*}
where $p:(K \sm L)^{(k)} \ra K^{\sm k}\sm L^{(k)}$ is the $\Sigma_k$--equivariant projection.
\end{enumerate}
\end{thm}

\begin{cor}  The natural transformation of functors of $X$,
$$ \Sinfty \eta: \Sinfty \MapT(L,X) \ra \Sinfty \MapT(K\sm L, K\sm X)$$
induces, on $k^{th}$ Taylor coefficients, the $\Sigma_k$--equivariant map of spectra
$$ \MapS(L^{(k)}, S) \xra{\eta} \MapS(K^{\sm k} \sm L^{(k)},   K^{\sm k}) \xra{p^*} \MapS((K \sm L)^{(k)},   K^{\sm k}).$$ 
\end{cor}

\begin{proof}[Proof of \thmref{smash theorem}] Given a surjection $\alpha: {\mathbf  n} \ra {\mathbf  m}$, we define 
$$ \eta_{\alpha}: \MapS(L^{\sm m}, X^{\sm n}) \ra \MapS((K \sm L)^{\sm m}, (K \sm X)^{\sm n})$$
to be the composite 
\begin{equation*}
\begin{split} 
\MapS(L^{\sm m}, X^{\sm n}) & \xra{\eta} \MapS(K^{\sm m} \sm L^{\sm m}, K^{\sm m} \sm X^{\sm n}) \\
& \xra{\alpha^*} \MapS(K^{\sm m} \sm L^{\sm m}, K^{\sm n} \sm X^{\sm n}).
\end{split}
\end{equation*}
Note that $\eta_{\alpha}$ is natural in all variables.

Given surjections ${\mathbf  n}\xra{\beta}{\mathbf  m}\xra{\alpha}{\mathbf  l}$, one easily verifies that there is a commutative diagram
\begin{equation}\label{smash diagram}
\xymatrix{
 \MapS(L^{\sm m}, X^{\sm n})    \ar[r]^-{\eta_{\beta}} \ar[d]^{\alpha_*} & \MapS((K \sm L)^{\sm m}, (K \sm X)^{\sm n})
 \ar[d]^{\alpha_*} \\
 \MapS(L^{\sm l}, X^{\sm n})    \ar[r]^-{\eta_{\alpha \circ \beta}} & \MapS((K \sm L)^{\sm l}, (K \sm X)^{\sm n})
  \\
 \MapS(L^{\sm l}, X^{\sm m})    \ar[r]^-{\eta_{\alpha}} \ar[u]_{\beta^*} & \MapS((K \sm L)^{\sm l}, (K \sm X)^{\sm m}).
 \ar[u]_{\beta^*}}
\end{equation}
Recall that $P^L_{\infty}(X)$ is defined as an equalizer.  A first consequence of (\ref{smash diagram}) is that there is a commutative diagram

\begin{equation*}
\xymatrix{\displaystyle \prod_{n} \MapS(L^{\sm n}, X^{\sm n})
 \ar[d]^{ \prod_n \eta} \ar[r] \ar@<2ex>[r]&
 \displaystyle \prod_{{\mathbf  m}\xra{\alpha}{\mathbf  l}} \MapS(L^{\sm l}, X^{\sm m}) \ar[d]^{ \prod_{\alpha} \eta_{\alpha}} \\
\displaystyle \prod_{n} \MapS((K \sm L)^{\sm n}, (K \sm X)^{\sm n}) \ar[r] \ar@<2ex>[r] &
\displaystyle \prod_{{\mathbf  m}\xra{\alpha}{\mathbf  l}} \MapS((K \sm L)^{\sm l}, (K \sm X)^{\sm m}). 
}
\end{equation*}
By taking equalizers, this then induces a filtration preserving natural map  $$\eta:P^L_{\infty}(X) \ra P^{K \sm L}_{\infty}(K \sm X),$$
and thus a map on the associated towers.

To identify the induced map on fibers, we first note that 
$$ \eta_{\alpha}: \MapS(L^{\sm m}, X^{\sm n}) \ra \MapS((K \sm L)^{\sm m}, (K \sm X)^{\sm n})$$
also equals the composite 
\begin{equation*}
\begin{split} 
\MapS(L^{\sm m}, X^{\sm n}) & \xra{\eta} \MapS(K^{\sm n} \sm L^{\sm m}, K^{\sm n} \sm X^{\sm n}) \\
& \xra{\alpha_*} \MapS(K^{\sm m} \sm L^{\sm m}, K^{\sm n} \sm X^{\sm n}).
\end{split}
\end{equation*}
This makes it clear that the $\Sigma_k$--equivariant map
$$ \lim_{\alpha \in (\mathbf k \da \E_{k-1})}\eta_{\alpha}: \MapS(\Delta_k(L), X^{\sm k}) \ra 
\MapS(\Delta_k(K \sm L), (K \sm X)^{\sm k})$$
can be identified as the composite 
\begin{equation*}
\begin{split} 
\Map_{\Sp}(\Delta_k(L), X^{\sm k}) & \xra{\eta} \Map_{\Sp}(K^{\sm k} \sm \Delta_k(L), K^{\sm k} \sm X^{\sm k}) \\
& \ra \Map_{\Sp}(\Delta_k(K \sm L), (K \sm X)^{\sm k}),
\end{split}
\end{equation*}
and the last part of the theorem follows.

Finally, statement (1) of the theorem is a consequence of the following commutative diagram:
\begin{equation*}
\xymatrix{
\Sinfty \MapT(L,X) \ar[d]^{\Sinfty \xi_k} \ar[r]^-{\eta} &
\Sinfty \MapT(K \sm L, K \sm X) \ar[d]^{\Sinfty \xi_k}  \\
\Sinfty \MapT(L^{\sm k},X^{\sm k}) \ar[d]^s \ar[r]^-{\eta} &
\Sinfty \MapT((K \sm L)^{\sm k}, (K \sm X)^{\sm k}) \ar[d]^s  \\
 \MapS(L^{\sm k},X^{\sm k})  \ar[r]^-{\eta} &
 \MapS((K \sm L)^{\sm k}, (K \sm X)^{\sm k}).  \\
}
\end{equation*}
\end{proof}

\subsection{The evaluation theorem}  The next of our general theorems concerns the compatiblity of the Arone tower with generalized evaluation maps. 

We have maps
$$ \epsilon: K \sm \MapT(K \sm L,X) \ra \MapT(L, X)$$
defined as the adjoint to the evaluation map
$$ \epsilon_{\T}: K \sm L \sm \MapT(K \sm L,X) \ra X$$
already discussed in \secref{smash products of spectra}.
  
Using the stable evaluation $\epsilon_{\Sp}$, we similarly have maps
$$ \epsilon: K \sm \MapS(K \sm L,X) \ra \MapS(L, X).$$
These unstable and stable generalized evaluation maps are compatible under stablization, courtesy of \lemref{eval lemma}.

\begin{thm} \label{eval theorem}  There are natural maps  of towers
$$\epsilon:K \sm \Tilde{P}^{K \sm L}(X) \ra \Tilde{P}^L(X)$$
 with the following properties. 
\begin{enumerate}
\item There is a commutative diagram in $\Sp$:
$$
\xymatrix{
\Sinfty K \sm \MapT(K \sm L,X)    \ar[d]^{\epsilon} \ar[rr]^-{1_K \sm \Tilde{e}^{K \sm L}(X)} & &
K \sm \Tilde{P}^{K \sm L}(X)   \ar[d]^{\epsilon}  \\
\Sinfty\MapT(L,X) \ar[rr]^-{\Tilde{e}^{L}(X)} & &
\Tilde{P}^{L}( X). 
}
$$
\item The induced map on $k^{th}$ fibers  is the  
composite  
\begin{equation*}
\begin{split} K \sm \Map_{\Sp}^{\Sigma_k}((K \sm L)^{( k)}, X^{\sm k}) & \xra{d^*} K \sm \Map_{\Sp}^{\Sigma_k}(K \sm L^{( k)}, X^{\sm k}) \\
& \xra{\epsilon} \Map_{\Sp}^{\Sigma_k}(L^{( k)}, X^{\sm k}),
\end{split}
\end{equation*}
where the first map is induced by the $\Sigma_k$--equivariant map of spaces
$$ d: K \sm L^{(k)} \ra (K \sm L)^{(k)}$$
which arises by embedding $K$ diagonally in $K^{\sm k}$.
\end{enumerate}
\end{thm}

\begin{cor}  The natural transformation of functors of $X$,
$$ \Sinfty \epsilon: \Sinfty K \sm \MapT(K \sm L,X) \ra \Sinfty \MapT(L, X)$$
induces, on $k^{th}$ Taylor coefficients, the $\Sigma_k$--equivariant map of spectra
$$ K \sm \MapS((K \sm L)^{(k)}, S) \xra{d^*} K \sm \MapS(K \sm L^{(k)}, S) \xra{\epsilon} \MapS(L^{(k)}, S).$$ 
\end{cor}

\begin{proof}[Proof of \thmref{eval theorem}] Given a surjection $\alpha: {\mathbf  n} \ra {\mathbf  m}$, we define 
$$ \epsilon_{\alpha}: K \sm \MapS((K \sm L)^{\sm m}, X^{\sm n}) \ra \MapS(L^{\sm m}, X^{\sm n})$$
to be the composite  
\begin{equation*}
\begin{split}
 K \sm \Map_{\Sp}(K^{\sm m} \sm L^{\sm m}, X^{\sm n}) & \xra{d^*} K \sm \Map_{\Sp}(K \sm L^{\sm m}, X^{\sm n}) \\
\xra{\epsilon} \Map_{\Sp}(L^{\sm m}, X^{\sm n}),
\end{split}
\end{equation*}
where the first map is induced by the diagonal $ d: K  \ra K^{\sm m}$. Note that $\epsilon_{\alpha}$ is natural in all variables.

Given surjections ${\mathbf  n}\xra{\beta}{\mathbf  m}\xra{\alpha}{\mathbf  l}$, one easily verifies that there is a commutative diagram
\begin{equation*}\label{evaluation diagram}
\xymatrix{
 K \sm \MapS((K \sm L)^{\sm m}, X^{\sm n})    \ar[r]^-{\epsilon_{\beta}} \ar[d]^{\alpha_*} & \MapS(L^{\sm m}, X^{\sm n}) \ar[d]^{\alpha_*} \\
 K \sm \MapS((K \sm L)^{\sm l}, X^{\sm n})    \ar[r]^-{\epsilon_{\alpha \circ \beta}} & \MapS(L^{\sm l}, X^{\sm n})  \\
 K \sm \MapS((K \sm L)^{\sm l}, X^{\sm m})    \ar[r]^-{\epsilon_{\alpha}} \ar[u]_{\beta^*} & \MapS(L^{\sm l}, X^{\sm m}).
 \ar[u]_{\beta^*}}
\end{equation*}
As in the proof of \thmref{smash theorem},
there is thus an induced natural map of towers $$\epsilon:K \sm \Tilde{P}^{K \sm L}(X) \ra \Tilde{P}^L(X).$$
The induced map on fibers is easy to identify. Note that 
the $\Sigma_k$--equivariant map
$$ \lim_{\alpha \in (\mathbf k \da \E_{k-1})}\epsilon_{\alpha}: K \sm \MapS(\Delta_k(K \sm L), X^{\sm k}) \ra 
\MapS(\Delta_k(L), X^{\sm k})$$
can be identified as the composite 
\begin{equation*}
\begin{split}
K \sm \Map_{\Sp}(\Delta_k(K \sm L), X^{\sm k}) & \ra K \sm \Map_{\Sp}(K \sm \Delta_k(L), X^{\sm k}) \\
&  \xra{\epsilon} \Map_{\Sp}(\Delta_k(L), X^{\sm k}),
\end{split}
\end{equation*}
where the first map is induced by the equivariant inclusion 
$$K \sm \Delta_k(L) \ra \Delta_k(K \sm L).$$
The last part of the theorem follows.

Statement (1) of the theorem follows by juxtaposing the two commutative diagrams

\begin{equation*}
\xymatrix{ 
\Sinfty K \sm \MapT(K \sm L,X) \ar[r]^-{1 \sm \xi_k} \ar[dr]^{\Delta \sm \xi_k} \ar[dd]^-{\epsilon} & \Sinfty K \sm \MapT((K \sm L)^{\sm k},X^{\sm k}) \ar[d]^-{\Delta \sm 1}
  \\
 & \Sinfty K^{\sm k} \sm \MapT((K \sm L)^{\sm k},X^{\sm k}) \ar[d]^-{\epsilon} \\
 \Sinfty \MapT(L, X) \ar[r]^-{\xi_k} & \Sinfty \MapT((K \sm L)^{\sm k}, (K \sm X)^{\sm k}),   
}
\end{equation*}
and 
\begin{equation*}
\xymatrix{ 
 \Sinfty K \sm \MapT((K \sm L)^{\sm k},X^{\sm k}) \ar[r]^-{1 \sm s} \ar[d]^-{\Delta \sm 1}& K \sm \MapS((K \sm L)^{\sm k},X^{\sm k})  \ar[d]^-{\Delta \sm 1}
  \\
\Sinfty K^{\sm k} \sm \MapT((K \sm L)^{\sm k},X^{\sm k}) \ar[d]^-{\epsilon} \ar[r]^-{1 \sm s} &  K^{\sm k} \sm \MapS((K \sm L)^{\sm k},X^{\sm k})  \ar[d]^-{\epsilon}  \\
 \Sinfty \MapT((K \sm L)^{\sm k}, (K \sm X)^{\sm k}) \ar[r]^-s   &
 \MapS((K \sm L)^{\sm k}, (K \sm X)^{\sm k}).  
}
\end{equation*}
\end{proof}

\section{The Product and Diagonal Theorems} \label{product section}

The goal of this section is to prove a general result about how the unstable homeomorphisms
$$ \MapT(K \vee L, X) = \MapT(K,X) \times \MapT(L,X)$$
lead to pairings among the associated Arone towers. One thus gets pairings  among the  spectral sequences which arise after applying any multiplicative cohomology theory to the towers.

In discussing and proving our result, it seems necessary to first say a little about the general theory of the smash product of towers.  The results here are unsurprising and presumably known, but we have been unable to find them explicitly in the literature.  

\subsection{Homotopy limits of spectra}
If $\J$ is a small category, $E\J_+: \J \ra \T$ is defined by letting $$E\J(j)_+ = B(\J \da j)_+.$$  Given a functor $Y: \J \ra \Sp$, its homotopy limit is the spectrum defined by the formula
$$ \holim_{\J}Y = \Map_{\Sp}^{\J}(E\J_+, Y).$$ 
This construction is natural with respect to both the functor and the category.  

We will make use of the following construction.  Given $X \in \Sp$ and $Y:\J \ra \Sp$ there is a natural map
\begin{equation} \label{smashlimitmap1} 
X \sm \holim_{\J} Y \ra \holim_{\J} (X \sm Y)
\end{equation}
defined as the adjoint to the composite
$$ \Map_{\Sp}^{\J}(E\J_+, Y) \ra \Map_{\Sp}^{\J}(X \sm E\J_+, X \sm Y) = \MapS( X, \Map_{\Sp}^{\J}(E\J_+, X \sm Y)).$$
Similarly, there is a natural map
\begin{equation} \label{smashlimitmap2} 
(\holim_{\J} Y) \sm X \ra \holim_{\J} (Y \sm X)
\end{equation}
The notation $X \times_Z Y$ will denote the homotopy pullback \ 
$ \holim \left( \def\objectstyle{\scriptstyle} \vcenter{
\xymatrix @-1.2pc {
 &
Y \ar[d]  \\
X \ar[r] &
Z
}}
\right).
$
The homotopy fiber of a map $Y \ra Z$ is then defined as $* \times_Z Y$.

\subsection{Smash products of towers of spectra} \label{towers of spectra}

If we let $\N$ denote the poset $0 < 1 < 2 < \dots$, then a tower of spectra can be regarded as a functor $P: \N^{op} \ra \Sp$.  Given towers $P$ and $Q$, respectively taking values in $\Sp \U$ and $\Sp \U^{\prime}$, the external smash product yields a bi--tower $P \sm Q: (\N \times \N)^{op} \ra \Sp(\U \times \U^{\prime})$.  When the context is clear, we will repress explicit notation for universes, and thus just write $P \sm Q: (\N \times \N)^{op} \ra \Sp$.

Now suppose given a general bi--tower $C: (\N \times \N)^{op} \ra \Sp$.  One then has an associated tower (which we still call $C$) with $k^{th}$ term given by 
$$ C_k = \holim_{i+j \leq k} C_{i,j}.$$
This tower can be understood homotopically via the next lemma.

\begin{lem} \label{bitower lemma} Given a bi--tower $C: (\N \times \N)^{op} \ra \Sp$, the diagram
\begin{equation*}
\xymatrix{
C_k \ar[d] \ar[r] &
\prod_{i+j=k} C_{i,j} \ar[d]  \\
C_{k-1} \ar[r] &
\prod_{i+j=k} C_{i-1,j} \times_{C_{i-1,j-1}} C_{i,j-1}
}
\end{equation*}
naturally homotopy commutes, and induces a natural homotopy equivalence on homotopy fibers. \\

\end{lem}

\begin{proof}[Sketch proof]
Let $\N_i \subset \N$ be the poset $0 < 1 < \dots < i$, and let $(\N \times \N)_k \subset \N \times \N$ be the poset $\{ (i,j) \ | \ i+j \leq k \}$.  There is a pushout of posets:
\begin{equation*} 
\xymatrix{
 \coprod_{i+j = k} \N_i \times \N_j - \{(i,j)\} \ar[d] \ar[r] &
(\N \times \N)_{k-1} \ar[d]  \\
 \coprod_{i+j = k} \N_i \times \N_j \ar[r] &
(\N \times \N)_{k}.
}
\end{equation*}
This induces a pullback of spectra
\begin{equation*}
\xymatrix{
C_k \ar[d] \ar[r] &
 \prod_{i+j = k} \holim_{\N_i \times \N_j} C \ar[d]  \\
C_{k-1} \ar[r] &
 \prod_{i+j = k} \holim_{\N_i \times \N_j - \{(i,j)\} } C
}
\end{equation*}
in which the two vertical maps are fibrations.  Now note that 
$$ \holim_{\N_i \times \N_j} C \ra C(i,j)$$
is an equivalence, as $(i,j)$ is the terminal object in $\N_i \times \N_j$, 
and 
$$ \holim_{\N_i \times \N_j - (i,j)} C \ra C_{i-1,j} \times_{C_{i-1,j-1}} C_{i,j-1}$$
is an equivalence, as $ \left\{ \def\objectstyle{\scriptstyle} \vcenter{
\xymatrix @-1.2pc {
(i-1,j) &
  \\
(i-1,j-1) \ar[r] \ar[u]&
(i,j-1)
}}
\right\}
$
is cofinal in $\N_i \times \N_j - \{(i,j)\}$. 
\end{proof}

As in the introduction, if $P$ is a tower, we let $F_k(P)$ denote the homotopy fiber of $P_k \ra P_{k-1}$.  Similarly, if $C$ is a bi--tower, we let $F_{i,j}(C)$ denote the homotopy fiber of $C_{i,j} \ra C_{i-1,j} \times_{C_{i-1,j-1}} C_{i,j-1}$.   This is the same thing as the iterated homotopy fiber of the square
\begin{equation*}
\xymatrix{
C_{i,j} \ar[d] \ar[r] &
C_{i,j-1} \ar[d]  \\
C_{i-1,j} \ar[r] &
C_{i-1,j-1}. 
}         
\end{equation*}
(See \cite[1.1b]{goodwillie2}.)
Specializing to the case when $C = P \sm Q$, use of construction (\ref{smashlimitmap2}), and then (\ref{smashlimitmap1}), thus yields the composite
\begin{equation} \label{fiber map} F_i(P) \sm F_j(Q) \ra F_i(P \sm F_j(Q)) \ra F_{i,j}(P\sm Q).
\end{equation}

\begin{lem}  $ F_i(P) \sm F_j(Q) \ra F_{i,j}(P\sm Q)$ is a natural weak equivalence of spectra. \\

\end{lem}
\begin{proof}[Sketch proof] As a fibration sequence in spectra is homotopy equivalent to a cofibration sequence, smashing a fibration sequence with a spectrum yields a sequence equivalent to a fibration sequence.  Applying this principle to each of the two maps in (\ref{fiber map}) shows each to be an equivalence.
\end{proof}

The two lemmas together have the following corollary. \\

\begin{cor} \label{fibprodcor} If $P$ and $Q$ are towers of spectra, the natural maps
$$ F_k(P\sm Q) \ra \prod_{i+j=k} F_{i,j}(P\sm Q) \la \prod_{i+j=k} F_i(P)\sm F_j(Q)$$
are weak equivalences.
\end{cor}

We now turn to a discussion of spectral sequences.  Suppose $h$ is a spectrum.  Applying $h_*$ and $h^*$ to a tower $P$ yields left half plane homology and cohomology spectral sequences $\{ E^r_{*,*}(P)\}$ and $\{ E_r^{*,*}(P)\}$ with 
$$ E^1_{-k,*}(P) = h_{*-k}(F_k(P))$$
and 
$$ E_1^{-k,*}(P) = h^{*-k}(F_k(P)).$$

Now suppose $h$ is a ring spectrum.  Then the last corollary implies that, given towers $P$ and $Q$,  there are natural pairings
$$ E^1_{-i,m}(P)\otimes E^1_{-j,n}(Q)\ra E^1_{-(i+j),m+n}(P\sm Q)$$
and 
$$ E_1^{-i,m}(P)\otimes E_1^{-j,n}(Q)\ra E_1^{-(i+j),m+n}(P\sm Q).$$
One then has:
\begin{thm}  The pairings extend to pairings of spectral sequences.
\end{thm}
The pairings have the expected properties.  For example, in the cohomology spectral sequence, the pairings are related to the pairing of filtered groups $$\colim_i h^*(P_i) \otimes \colim_j h^*(Q_j)  \ra \colim_k h^*((P\sm Q)_k)$$
at the level of $E_{\infty}$.

A careful proof of this theorem seems to not appear in the literature.  However, a proof can be constructed in a straightforward manner by carefully mimicking the discussion on pages 660--668 of G.W.Whitehead's book \cite{whitehead} where he discusses pairings in spectral sequences associated to products of filtered spaces.  The translation into our setting involves arguments of the sort given in the proofs of our two lemmas, but no other new ideas.

\subsection{The Product Theorem}

The homeomorphism 
$$ \MapT(K \vee L, X)_+ =  \MapT(K,X)_+ \sm \MapT(L,X)_+   $$
suggests that there should be a compatible weak equivalence between the two towers $P^{K \vee L}(X)$ and $P^K(X)\sm P^L(X)$.  We will shortly see that this is the case.

To be computationally useful, we will need to also identify the induced weak equivalence on fibers.  By \corref{fibprodcor}, we have homotopy equivalences
\begin{equation*} \label{iso1} 
\begin{split}
F_k(P^K(X)\sm P^L(X)) & \xra{\sim} \prod_{i+j=k} F_{i,j}(P^K(X)\sm P^L(X)) \\
& \xla{\sim} \prod_{i+j=k} \Map_{\Sp}^{\Sigma_i}(K^{(i)}, X^{\sm i})\sm \Map_{\Sp}^{\Sigma_j}(L^{(i)}, X^{\sm j}).
\end{split}
\end{equation*}
Meanwhile, the $\Sigma_k$--equivariant homeomorphism
\begin{equation*} \label{wedge1} \bigvee_{i+j=k} {\Sigma_k}_+ \sm_{\Sigma_i \times \Sigma_j} (K^{(i)} \sm L^{(j)}) = (K \vee L)^{(k)}
\end{equation*}
induces an isomorphism
\begin{equation*} \label{iso2} F_k^{K\vee L}(X) = \prod_{i+j=k} \Map_{\Sp}^{\Sigma_i \times \Sigma_j}(K^{(i)}\sm L^{(j)},  X^{\sm k}).
\end{equation*}
To state our main product theorem it is convenient to define a bi--tower \\  $P^{K,L}(X)$ by the formula
$$ P_{i,j}^{K,L}(X) = \Map_{\Sp}^{\E_i \times \E_j}(K^{\sm} \sm L^{\sm}, \add^*X^{\sm}).$$
Here $\add^*X^{\sm}: (\E \times \E)^{op} \ra \T$ is defined by $\add^*X^{\sm}(\mathbf i, \mathbf j) = X^{\sm (i+j)}$.

\begin{thm} \label{KLprodthm}  There are natural homotopy equivalences of towers 
$$ P^{K\vee L}(X) \xra{\sim} P^{K,L}(X) \xla{\sim} P^K(X)\sm P^L(X)$$
with the following properties.
\begin{enumerate}
\item There is a commutative diagram in $\Sp$:
{\small
$$
\xymatrix{
\Sinfty \MapT(K \vee L,X)_+    \ar@{=}[dd] \ar[rr]^-{e^{K \vee L}(X)} & &
P^{K \vee L}(X)   \ar[d]^{\wr}  \\
&& P^{K,L}(X) \\
\Sinfty \MapT(K,X)_+ \sm \MapT(L,X)_+ \ar[rr]^-{e^K(X) \sm e^{L}(X)} & &
P^K(X) \sm P^{L}( X) \ar[u]_{\wr}. 
}
$$}

\item  The induced equivalences on $k^{th}$ fibers
$$ F_k^{K\vee L}(X) \xra{\sim} F_k(P^{K,L}(X)) \xla{\sim} F_k(P^K(X) \sm P^L(X))$$
\end{enumerate}
fits into a commutative diagram for each $i+j=k$:
{\footnotesize
\begin{equation*}
\xymatrix{
F_k^{K\vee L}(X) \ar[d] \ar[r]^-{\sim} & F_k(P^{K,L}(X)) \ar[d] & F_k(P^K(X) \sm P^L(X)) \ar[l]_-{\sim} \ar[d] \\
\Map_{\Sp}^{\Sigma_i \times \Sigma_j}(K^{(i)}\sm L^{(j)},  X^{\sm k}) \ar@{=}[d] \ar[r]^-{\sim} & F_{i,j}(P^{K,L}(X)) &  F_{i,j}(P^{K}(X)\sm P^L(X)) \ar[l]_-{\sim} \\
\Map_{\Sp}^{\Sigma_i \times \Sigma_j}(K^{(i)}\sm L^{(j)},  X^{\sm k}) & & *\txt{ $\Map_{\Sp}^{\Sigma_i}(K^{(i)},  X^{\sm i})$ \\ $\sm \Map_{\Sp}^{\Sigma_j}(L^{(j)},  X^{\sm j})$,} \ar[ll]_-{\sim} \ar[u]_-{\wr}
}
\end{equation*}}
where the bottom map is the evident smash product map. \\
\end{thm}

We will use the notation $\mu: P^{K \vee L}(X) \xra{\sim} P^K(X) \sm P^L(X)$ to denote the weak natural equivalence of the theorem.

\begin{rem}  In stating this theorem, we have continued to repress notation for universes. However, we hope it is understood that, if $P^K(X)$ is a tower in $\Sp \U$, and $P^L(X)$ is a tower in $\Sp \U^{\prime}$, then the maps and objects in the theorem are living in $\Sp (\U \oplus \U^{\prime})$. \\
\end{rem}

The key to our product theorem is the following observation, hinted at in (\ref{wedge1}) above.
Let $\add: \E \times \E \ra \E$ be the functor defined by $\add(\mathbf i, \mathbf j) = (\mathbf{i+j})$.  This induces a functor
$$  \add^*: \T^{\E} \ra \T^{\E \times \E}$$
by pullback.  This has a left adjoint $\add_*: \T^{\E \times \E} \ra \T^{\E }$
explicitly given by 
$$\add_*X(\mathbf k) = \colim_{\mathbf k \da \E \times \E} X$$
where $\mathbf k \da \E \times \E$ is the category with objects all triples $(\mathbf i, \mathbf j, \alpha)$ where \\
$\alpha:\mathbf k \ra \mathbf{ (i+j)}$ is a surjection.  Then one has

\begin{lem} The inclusions $K^{\sm m} \sm L^{\sm n} \subset (K \vee L)^{\sm m+n}$ induce an isomorphism in $\T^{\E}$,
$$\add_*(K^{\sm}\sm L^{\sm}) = (K \vee L)^{\sm},$$
 and this restricts to give isomorphisms for all $k$
$$ \add_*(\colim_{i+j \leq k} (K_i^{\sm}\sm L_j^{\sm})) = (K \vee L)_k^{\sm}.$$
\end{lem}

We also record the following fact.

\begin{lem} $\displaystyle \hocolim_{i+j \leq k} (K_i^{\sm}\sm L_j^{\sm}) \ra \colim_{i+j \leq k} (K_i^{\sm}\sm L_j^{\sm})$ is an equivalence in $\T^{\E \times \E}$. 
\end{lem}
\begin{proof} Let $(\N \times \N)_k$ denote the full subcategory of $\N \times \N$ with objects all $(i,j)$ such that $i+j \leq k$.  It is easy to see that, for all $n$, $K$, and $L$, the functor on $\N \times \N$ sending $(i,j)$ to the space $K_i^{\sm}(n) \sm L_j^{\sm}(n)$ is a cofibrant object in the model category structure on $\T^{(\N \times \N)_k}$ described in  \cite[\S 10.13]{dwyerspalinski}. The lemma follows.
\end{proof}

The maps in the theorem are now easy to define.  A natural equivalence 
$$ P_k^{K\vee L}(X) \ra P_k^{K,L}(X)$$
is defined by
\begin{equation*}
\begin{split}
P_k^{K\vee L}(X) & = \Map_{\Sp}^{\E}((K\vee L)^{\sm}_k, X^{\sm}) \\
  & = \Map_{\Sp}^{\E}(\add_*(\colim_{i+j \leq k} (K_i^{\sm}\sm L_j^{\sm})), X^{\sm})\\
 & = \Map_{\Sp}^{\E \times \E}(\colim_{i+j \leq k} (K_i^{\sm}\sm L_j^{\sm}), \add^*X^{\sm})\\
 & \xra{\sim} \Map_{\Sp}^{\E \times \E}(\hocolim_{i+j \leq k} (K_i^{\sm}\sm L_j^{\sm}), \add^*X^{\sm}) \\
 & = \holim_{i+j \leq k} \Map_{\Sp}^{\E \times \E}(K_i^{\sm}\sm L_j^{\sm}, \add^*X^{\sm}) \\
 & = \holim_{i+j \leq k} P_{i,j}^{K,L}(X) = P_k^{K,L}(X).
\end{split}
\end{equation*}
The map of towers  \ $ P^K(X)\sm P^L(X) \ra P^{K,L}(X)$ \ is even more evident.  It is the map on towers induced by the map of bi--towers
$$ \sm: \Map_{\Sp}^{\E_i}(K^{\sm}, X^{\sm})\sm \Map_{\Sp}^{\E_j}(L^{\sm}, X^{\sm}) \ra \Map_{\Sp}^{\E_i \times \E_j}(K^{\sm} \sm L^{\sm}, \add^*X^{\sm}).$$
(This is not so evidently a homotopy equivalence, but, we will learn that it is.) 

The theorem is now easily proved.

First, we check that the diagram in (1) commutes.  This will follow if we  verify that, for all $i+j = k$, the diagram
{\footnotesize 
$$
\xymatrix{
\Sinfty \MapT(K \vee L,X)_+    \ar@{=}[dd] \ar[rr]^-{e^{K \vee L}(k)} & &
\MapS((K \vee L)^{\sm k}, X^{\sm k})   \ar[d]  \\
&& \MapS(K^{\sm i} \sm  L^{\sm j}, X^{\sm k}) \\
\Sinfty \MapT(K,X)_+ \sm \MapT(L,X)_+ \ar[rr]^-{e^K(i) \sm e^{L}(j)} & &
\MapS(K^{\sm i},X^{\sm i}) \sm \MapS(L^{\sm j},X^{\sm j})  \ar[u]. 
}
$$}
commutes, where we have written $e^K(i)$ for $e^K(X,i)$, etc.   But this diagram commutes, as one easily checks that the diagram of spaces 
{\small 
$$
\xymatrix{
 \MapT(K \vee L,X)    \ar@{=}[dd] \ar[rr]^-{\xi_k} & &
\MapT((K \vee L)^{\sm k}, X^{\sm k})   \ar[d]  \\
&& \MapT(K^{\sm i} \sm  L^{\sm j}, X^{\sm k}) \\
\MapT(K,X) \times  \MapT(L,X) \ar[rr]^-{\xi_i \times \xi_j} & &
 \MapT(K^{\sm i},X^{\sm i}) \times \MapT(L^{\sm j},X^{\sm j})  \ar[u]. 
}
$$}
commutes, and then \lemref{smash lemma} implies that the diagram 

{\small 
$$
\xymatrix{
\Sinfty \MapT((K \vee L)^{\sm k}, X^{\sm k})    \ar[d] \ar[r]^-{s} &
\MapS((K \vee L)^{\sm k}, X^{\sm k})   \ar[d]  \\
\Sinfty \MapT(K^{\sm i} \sm  L^{\sm j}, X^{\sm k}) \ar[r]^-{s} & \MapS(K^{\sm i} \sm  L^{\sm j}, X^{\sm k}) \\
\Sinfty \MapT(K^{\sm i},X^{\sm i}) \sm \MapT(L^{\sm j},X^{\sm j}) \ar[r]^-{s \sm s} \ar[u]  &
\MapS(K^{\sm i},X^{\sm i}) \sm \MapS(L^{\sm j},X^{\sm j})  \ar[u]. 
}
$$}
commutes.

The lower rectangle in (2) commutes by inspection, the top left square commutes by definition, and the top right square commutes by naturality.   In this diagram, the two downward arrows are equivalences arising from \lemref{bitower lemma}, the top rightward arrow is an equivalence as it is induced by an equivalence of towers, and the lower leftward arrow is an equivalence by \corref{smash cor}.  It follows that all the other arrows here are equivalences, as asserted.  In particular, the map of towers $P^K(X) \sm P^L(X) \ra P^{K,L}(X)$ induces equivalences on all fibers, thus (inducting up the tower) is itself a homotopy equivalence.

\subsection{The Diagonal Theorem}

Let $\nabla: K \vee K \ra K$ be the fold map.  Since the diagonal on $\MapT(K,X)$ has a factorization
\begin{equation*}
\xymatrix{
\MapT(K,X) \ar[dr]^-{\nabla^*} \ar[rr]^-{\Delta} & &
\MapT(K,X) \times \MapT(K,X),  \\
& \MapT(K \vee K,X) \ar@{=}[ur] &
}
\end{equation*}
our Product Theorem has  consequences for the diagonal map.

Let $\Psi: P^K(X) \ra P^K(X) \sm P^K(X)$ be the weak natural transformation 
$$ P^K(X) \xra{\nabla^*} P^{K \vee K}(X) \xra{\mu} P^K(X) \sm P^K(X).$$

\begin{thm} The weak natural transformation $\Psi$ has the following properties.
\begin{enumerate}
\item There is a commutative diagram of weak natural transformations:
$$
\xymatrix{
\Sinfty \MapT(K,X)_+    \ar[d]^{\Delta} \ar[rr]^-{e^{K}(X)}  & &
P^{K}(X)   \ar[d]^{\Psi}  \\
\Sinfty \MapT(K,X)_+ \sm \MapT(K,X)_+ \ar[rr]^-{e^K(X) \sm e^{K}(X)} & &
P^K(X) \sm P^{K}( X). 
}
$$
\item Via the natural weak equivalences $F_k^K(X) \simeq \MapS(K^{(k)}, X^{\sm k})_{h\Sigma_k}$ and 
$$F_k(P^K(X) \sm P^K(X)) \simeq \prod_{i+j = k} \MapS(K^{(i)}, X^{\sm i})_{h\Sigma_i}\sm \MapS(K^{(j)}, X^{\sm j})_{h\Sigma_j},$$
the map induced by $\Psi$ on $k^{th}$ fibers is the product, over $i+j = k$, of the composites of the weak natural transformations
\begin{equation*}
\begin{split}
\MapS(K^{(k)}, X^{\sm k})_{h\Sigma_k} & \xra{Tr_{\Sigma_i \times \Sigma_j}^{\Sigma_k}}
\MapS(K^{(k)}, X^{\sm k})_{h(\Sigma_i \times \Sigma_j)} \\
 & \xra{\pi^*} \MapS(K^{(i)}\sm K^{(j)}, X^{\sm k})_{h(\Sigma_i \times \Sigma_j)} \\ 
& \xla{\sim} \MapS(K^{(i)}, X^{\sm i})_{h\Sigma_i}\sm \MapS(K^{(j)}, X^{\sm j})_{h\Sigma_j}, 
\end{split}
\end{equation*}
where $\pi:K^{(i)} \sm K^{(j)} \ra K^{(k)}$ is the projection.
\end{enumerate}

\end{thm}  

Property (1) follows immediately from property (1) of  \thmref{KLprodthm}.  To see that (2) follows from \thmref{KLprodthm}(2), we first observe that there is a factorization
$$ \xymatrix{
K^{(i)} \sm K^{(j)} \ar[rr]^{\pi} \ar[dr] & & K^{(k)}, \\
& (K \vee K)^{(k)} \ar[ur]^{\nabla^{(k)}}
}$$
and thus a commutative diagram

\begin{equation*}
\xymatrix{
\Map_{\Sp}^{\Sigma_k}(K^{(k)}, X^{\sm k}) \ar[d]^{\nabla^*} \ar[r] &
\Map_{\Sp}^{\Sigma_i \times \Sigma_j}(K^{(k)}, X^{\sm k}) \ar[d]^{\pi^*}  \\
\Map_{\Sp}^{\Sigma_k}((K \vee K)^{(k)}, X^{\sm k}) \ar[r] &
\Map_{\Sp}^{\Sigma_i \times \Sigma_j}(K^{(i)}\sm K^{(j)}, X^{\sm k}).
}
\end{equation*}
Thus \thmref{KLprodthm}(2) implies that the map on fibers can be identified with the product, over $i+j= k$, of the right vertical composites in the commutative diagrams of weak natural transformations
\begin{equation*}
\xymatrix{
\Map_{\Sp}(K^{(k)}, X^{\sm k})_{h\Sigma_k} \ar[d]^{Tr} \ar[r]^{\Phi}_{\sim} &
 \Map_{\Sp}^{\Sigma_k}(K^{(k)}, X^{\sm k})\ar[d]  \\
\Map_{\Sp}(K^{(k)}, X^{\sm k})_{h(\Sigma_i \times \Sigma_j)} \ar[r]^{\Phi}_{\sim} \ar[d]^{\pi^*} &
\Map_{\Sp}^{\Sigma_i \times \Sigma_j}(K^{(k)}, X^{\sm k}) \ar[d]^{\pi^*} \\
\Map_{\Sp}(K^{(i)}\sm K^{(j)}, X^{\sm k})_{h(\Sigma_i \times \Sigma_j)}  \ar[r]^{\Phi}_{\sim}  &
\Map_{\Sp}^{\Sigma_i \times \Sigma_j}(K^{(i)}\sm K^{(j)}, X^{\sm k}) \\
*\txt{$\MapS(K^{(i)}, X^{\sm i})_{h\Sigma_i}$ \\ $\sm \MapS(K^{(j)}, X^{\sm j})_{h\Sigma_j}$} \ar[r]^{\Phi \sm \Phi}_{\sim} \ar[u]_{\wr} &
*\txt{$\Map_{\Sp}^{\Sigma_i}(K^{(i)}, X^{\sm i})$ \\ $\sm \Map_{\Sp}^{\Sigma_j}(K^{(j)}, X^{\sm j}).$} \ar[u]_{\wr} \\
}
\end{equation*}

\section{Compatibility results} \label{compatibility section}

The following propositions and corollaries state that our various natural transformations are compatible in the expected ways.  The three propositions are easily verified by directly checking their definitions.

\begin{prop} \label{etaepsilonprop}  For all $J$, $K$, and $L$, the diagram of natural transformations of towers
\begin{equation*}
\xymatrix{
K \sm \Tilde{P}^{K \sm L}(X) \ar[d]^{\epsilon} \ar[r]^-{1 \sm \eta} &
K \sm \Tilde{P}^{J\sm K \sm L}(J \sm X)  \ar[r]^{\tau}_{\sim} & K \sm \Tilde{P}^{K\sm J \sm L}(J \sm X) \ar[d]^{\epsilon}  \\
\Tilde{P}^L(X) \ar[rr]^-{\eta} & & \Tilde{P}^{J \sm L}(J \sm X)
}
\end{equation*}
commutes, where $\tau$ is the isomorphism induced by the twist map $J \sm K \ra K \sm J$.
\end{prop}

\begin{prop} \label{etamuprop}  For all $J$, $K$, and $L$, the diagram of weak natural transformations of towers
\begin{equation*}
\xymatrix{
P^{K \vee L}(X) \ar[d]^{\mu}_{\wr} \ar[r]^-{\eta} &
P^{J\sm (K \vee L)}(J \sm X)  \ar@{=}[r] & P^{(J \sm K) \vee(J \sm L)}(J \sm X) \ar[d]^{\mu}_{\wr}  \\
P^K(X) \sm P^L(X) \ar[rr]^-{\eta \sm \eta} & & P^{J \sm K}(J \sm X) \sm P^{J \sm L}(J \sm X)
}
\end{equation*}
commutes.
\end{prop}

\begin{cor} \label{etapsicor}  For all $J$ and $K$, the diagram of weak natural transformations of towers
\begin{equation*}
\xymatrix{
P^{K }(X) \ar[d]^{\Psi} \ar[rr]^-{\eta} & &
P^{J\sm K}(J \sm X)  \ar[d]^{\Psi}  \\
P^K(X) \sm P^K(X) \ar[rr]^-{\eta \sm \eta} & & P^{J \sm K}(J \sm X) \sm P^{J \sm K}(J \sm X)
}
\end{equation*}
commutes.
\end{cor}

\begin{ex}  A consequence of this corollary is that the truth of \corref{Sn diagonal cor}, for all $n < \infty$, implies that \corref{Sn diagonal cor} is true when $n = \infty$. \\
\end{ex}

\begin{prop} \label{epsilonmuprop}  For all $J$, $K$, and $L$, the diagram of weak natural transformations of towers
\begin{equation*}
\xymatrix{
J \sm \Tilde{P}^{J \sm (K \vee L)}(X) \ar[d]^{\Delta \sm 1} \ar[rr]^-{\epsilon} & &
\Tilde{P}^{K \vee L}(X) \ar[dd]^{\mu}_{\wr}  \\
J \sm J \sm \Tilde{P}^{(J \sm K) \vee (J \sm L)} \ar[d]^{1 \sm \mu}_{\wr} & & \\
J \sm J \sm \Tilde{P}^{J \sm K} \sm P^{J \sm L}
\ar[r]^-{1 \sm \tau \sm 1}_-{\sim} & J \sm \Tilde{P}^{J \sm K} \sm J \sm \Tilde{P}^{J \sm L} \ar[r]^{\epsilon \sm \epsilon} & \Tilde{P}^K(X) \sm \Tilde{P}^L(X)
}
\end{equation*}
commutes, where $\Delta$ is the diagonal and $\tau$ is the twist isomorphism.
\end{prop}

\begin{cor} \label{epsilonpsicor}  For all $J$ and $K$, the diagram of weak natural transformations of towers
\begin{equation*}
\xymatrix{
J \sm \Tilde{P}^{J \sm K}(X) \ar[d]^{\Delta \sm \Psi} \ar[rr]^-{\epsilon} & &
\Tilde{P}^{K}(X) \ar[d]^{\Psi}  \\
J \sm J \sm \Tilde{P}^{J \sm K} \sm \Tilde{P}^{J \sm K}
\ar[r]^-{1 \sm \tau \sm 1}_-{\sim} & J \sm \Tilde{P}^{J \sm K} \sm J \sm \Tilde{P}^{J \sm K} \ar[r]^{\epsilon \sm \epsilon} & \Tilde{P}^K(X) \sm \Tilde{P}^K(X)
}
\end{equation*}
commutes, where $\Delta$ is the diagonal and $\tau$ is the twist isomorphism.
\end{cor}

\begin{ex}  With $J = S^1$, a typical consequence of this would be the following.  Let $\{E_r^{*,*}(K,X)\}$ be the spectral sequence obtained from $P^K(X)$ by applying a cohomology theory with products.  Then, because the reduced diagonal $\Delta: S^1 \ra S^1 \sm S^1$ is null, one deduces that 
$$ E_r(\epsilon): E_r^{*,*+1}(K,X) \ra E_r^{*,*}(\Sigma K,X)$$
is zero on the algebra decomposables in $E_r^{*,*}(K,X)$.  If the cohomology theory satisfies a Kunneth theorem (e.g.\ it is ordinary cohomology with field coefficients), the pinch map on $S^1$ shows that $E_r^{*,*}(\Sigma K,X)$ is a Hopf algebra, and formal manipulations then also imply that the image of $E_r(\epsilon)$ is contained in the primitives. \\
\end{ex}

\section{Little cubes and an explicit S--duality map} \label{duality section}

\subsection{Basic constructions with little cubes} \label{little cubes}

Let $I$ be the interval $[-1,1]$, and let $\C(n,1)$ be the space of `little $n$--cubes', the space of embeddings $I^n \ra I^n$ which are products of $n$ affine orientation preserving maps from $I$ to itself.  Then $\C(n,k)$ is defined to be the subspace of $\C(n,1)^k$ consisting of $k$--tuples of little $n$--cubes whose images have disjoint interiors.  Thus a point $c \in \C(n,k)$ can be viewed an embedding $c: \coprod_{i = 1}^k \stackrel{\circ}{I^n} \ra \stackrel{\circ}{I^n}$ of a special form.

Given a space $Z$, let $F(Z,k)\subseteq Z^k$ denote the configuration space of $k$ distinct points in $Z$:
 $$F(Z,k) = \{(z_1,\dots, z_k) \ | \ z_i \neq z_j \text{ if }i \neq j \}.$$  It is well known and easy to prove that the map 
$$ \C(n,k) \ra F(\stackrel{\circ}{I^n}, k),$$
sending a $k$--tuple of little $n$--cubes, 
$$(c_1,\dots,c_k),$$ to their centers, $$(c_1(\mathbf 0),\dots,c_k(\mathbf 0)),$$ is a $\Sigma_k$--equivariant homotopy equivalence.  

A point in $\C(n,k)$ provides a tubular neighborhood around the $0$-dimensional submanifold of $\stackrel{\circ}{I^n}$, consisting of the $k$ center points.  This suggests the following construction.  Given a point in $\C(n,k)$, \ 
$c: \coprod_{i = 1}^k \stackrel{\circ}{I^n} \longrightarrow \stackrel{\circ}{I^n}$, 
let $$c^*: S^n \longrightarrow \bigvee_{i=1}^k S^n$$  be the associated Thom--Pontryagin collapse map.  Then define
$$ \alpha(n,k): \C(n,k)_+ \ra \MapT(S^n, \bigvee_{i=1}^k S^n)$$
by $\alpha(n,k)(c) = c^*$.  Note that this is $\Sigma_k$--equivariant.

The maps $\alpha(n,k)$ are the starting points for two other families of maps. 

Let $\delta(n,1): \C(n,1)_+ \sm S^n \ra S^n$ be the adjoint of $\alpha(n,1)$.  Then notice that the subspace $\C(n,k)_+ \sm \Delta_k(S^n) \subset \C(n,k)_+ \sm S^{nk}$ is sent to the basepoint by the map 
$$\delta(n,1)^{\sm k}: \C(n,k)_+ \sm S^{nk} \subset \C(n,1)^k_+ \sm S^{nk} \ra S^{nk}.$$  Thus 
$\delta(n,1)^{\sm k}$ induces a $\Sigma_k$--equivariant map
$$\delta(n,k): \C(n,k)_+ \sm S^{n(k)} \ra S^{nk}.$$
A second family of $\Sigma_k$--equivariant maps
$$ \beta(m,n,k): S^m \sm \C(m+n,k)_+  \ra S^{mk} \sm \C(n,k)_+$$
is then defined by the following diagram:
\begin{equation*}
\xymatrix{
 S^m \sm \C(m+n,k)_+  \ar@{^{(}->}[d]^{\Delta \sm i} \ar[rr]^-{\beta(m,n,k)} &&
 S^{mk} \sm \C(n,k)_+  \ar@{^{(}->}[d]  \\
 S^{mk} \sm \C(m,1)^k_+ \sm \C(n,1)^k_+  \ar[rr]^-{(\delta(m,1)^{\sm k} \circ \tau) \sm 1} &&
 S^{mk} \sm \C(n,1)^k_+,
 }
 \end{equation*}
 where  $\tau: S^{mk} \sm \C(m,1)^k_+ \simeq \C(m,1)^k_+ \sm S^{mk}$ is the switch map, $\Delta: S^m \hra S^{mk}$ is the diagonal, and $i: \C(m+n,k) \hra \C(m,1)^k \times \C(n,1)^k$ is the map which regards each little $(m+n)$--cube as the product of a little $m$--cube with a little $n$--cube.

\subsection{The duality theorem and consequences}

The following duality theorem will be proved in \secref{duality proof}. \\

\begin{thm} \label{duality thm} The map $\delta(n,k): \C(n,k)_+ \sm S^{n(k)} \ra S^{nk}$ is an equivariant S--duality pairing.
\end{thm}

In other words, the stable adjoint
$$ \Tilde{\delta}(n,k): \Sinfty \C(n,k)_+ \ra \MapS(S^{n(k)}, S^{nk})$$
is a homotopy equivalence of $\Sigma_k$--spectra.

Given a space $X$, the map $\delta(n,k)$ induces a natural map
$$ \MapS(S^{nk},X^{\sm k}) \ra \MapS(\C(n,k)_+ \sm S^{n(k)}, X^{\sm k}).$$
A consequence of the duality theorem is that the adjoint of this,
\begin{equation} \label{duality iso} \C(n,k)_+ \sm \MapS(S^{nk}, X^{\sm k}) \ra \MapS(S^{n(k)}, X^{\sm k}),
\end{equation}
is a weak equivalence of $\Sigma_k$--spectra.  Passing to homotopy orbits, we have constructed the natural weak equivalence (\ref{Sn fiber1}) of the introduction:
$$ F^{S^n}_k(X) \simeq (\C(n,k)_+ \sm \MapS(S^{nk}, X^{\sm k}))_{h\Sigma_k}.$$
Using \cite[Thm.I.7.9 and Prop.VI.5.3]{lmms}, we then deduce (\ref{Sn fiber2}):
$$ F^{S^n}_k(X) \simeq  \C(n,k)_+ \sm_{\Sigma_k} (\Sigma^{-n}X)^{\sm k}.$$
Assuming the duality theorem, we now deduce the corollaries of the introduction from the corresponding theorems.

\begin{proof}[Proof of \corref{Sn smash cor}]
We specialize the Smashing Theorem to the case when $K = S^n$ and $L = S^m$.  We need to show that, under the equivalence (\ref{duality iso}), the description of the map on fibers given in the theorem, corresponds to the description given in the corollary.

By formal manipulation of adjunctions, we are asserting that there is a commutative diagram of $\Sigma_k$--spectra:
\begin{equation*}
\xymatrix{
\C(m,k)_+  \ar[rr]^-{\Tilde{\delta}(m,k)}_-{\sim}  \ar[dd]^-{i}  & & \MapS(S^{m(k)}, S^{mk}) \ar[d]^{\eta} \\
& & \MapS(S^{m(k)}\sm S^{nk}, S^{mk} \sm S^{nk}) \ar[d]^{p^*} \\ 
\C(m+n,k)_+ \ar[rr]^-{\Tilde{\delta}(m+n,k)}_-{\sim} & & \MapS((S^{m+n})^{(k)}, (S^{m+n})^{\sm k}).  
}
\end{equation*}
Here $i$ is the inclusion induced by multiplying all little $m$--cubes by the identity $n$--cube.

Again adjointing, we just need to check that there is a commutative diagram of $\Sigma_k$--spaces:
\begin{equation*}
\xymatrix{
\C(m,k)_+ \sm (S^{m+n})^{(k)} \ar[dd]^{1 \sm p} \ar[rr]^-{i \sm 1}  & &
\C(m+n,k)_+ \sm (S^{m+n})^{(k)}. \ar[dd]^{\delta(m+n,k)} \\
& &   \\
\C(m,k)_+ \sm S^{m(k)} \sm S^{nk} \ar[rr]^-{\delta(m,k)\sm 1} & & S^{mk} \sm S^{nk}=(S^{m+n})^{\sm k} 
}
\end{equation*}
This diagram, in turn, is a quotient of smash products of the diagram in the case when $k=1$.  This is the diagram
\begin{equation*}
\xymatrix{
\C(m,1)_+ \sm S^{m+n} \ar[dr]^-i \ar[r]^-{\delta(m,1) \sm 1} &  S^{m+n} \\
& \C(m+n)_+, \ar[u]_{\delta(m+n,1)} 
}
\end{equation*}
which is easily verified to be commutative.
\end{proof}

\begin{proof}[Proof of \corref{Sn diagonal cor}]
The corollary follows from the Diagonal Theorem, once we verify that, for all $i+j = k$, there is a 
commutative diagram of $\Sigma_i \times \Sigma_j$--spectra:
\begin{equation*}
\xymatrix{
\C(n,k)_+  \ar[rr]^-{\Tilde{\delta}(n,k)}_-{\sim} \ar[dd]^i &&
\MapS(S^{n(k)}, S^{nk})  \ar[d]^-{\pi^*} \\
& &  \MapS(S^{n(i)} \sm S^{n(j)}, S^{nk}) \\
(\C(n,i) \times \C(n,j))_+ \ar[rr]^-{\Tilde{\delta}(n,i) \sm \Tilde{\delta}(n,i)}_-{\sim}
& & \MapS(S^{n(i)}, S^{ni}) \sm  \MapS(S^{n(j)}, S^{nj})  \ar[u]_-{\sim} 
 }
\end{equation*}
Here $i$ is the inclusion which sends a $k$--tuple of little cubes to the first 
$i$ cubes and the last $j$ cubes.

This diagram of spectra commutes because there is a commutative diagram of $\Sigma_i \times \Sigma_j$--spaces:
\begin{equation*}
\xymatrix{
\C(n,k)_+ \sm S^{n(i)} \sm S^{m(j)} \ar[dd]^{1 \sm \pi} \ar[rr]^-{i \sm 1}  & &
\C(n,i)_+ \sm \C(n,j)_+ \sm S^{n(i)} \sm S^{n(j)} \ar[d]^{1 \sm \tau \sm 1} \\
& &
\C(n,i)_+ \sm S^{n(i)} \sm \C(n,j)_+ \sm S^{n(j)} \ar[d]^{\delta(n,i)\sm \delta(n,j)}   \\ 
\C(n,k)_+ \sm S^{n(k)}  \ar[rr]^-{\delta(n,k)} & & S^{nk} = S^{ni}\sm S^{nj}.
}
\end{equation*}
Here $\tau$ is the twist map.
\end{proof}

\begin{proof}[Proof of \corref{Sn evaluation cor}]
The corollary follows from the Evaluation Theorem, once we verify that there is a commutative diagram of $\Sigma_k$--spectra:
\begin{equation*}
\xymatrix{
S^m \sm \C(m+n,k)_+  \ar[rr]^-{1 \sm \Tilde{\delta}(m+n,k)}_-{\sim} \ar[dd]^{\beta(m,n,k)} & & 
S^m \sm \MapS((S^{m+n})^{(k)}, (S^{m+n})^{\sm k})  \ar[d]^-{1 \sm d^*} \\
& &  S^m \sm \MapS(S^m \sm S^{n(k)} , (S^{m+n})^{\sm k}) \ar[d]^-{\epsilon} \\
S^{mk} \sm \C(n,j))_+, \ar[rr]^-{\sim} &&  \MapS(S^{n(k)}, S^{mk} \sm S^{nk})
}
\end{equation*}
where the unlabelled horizontal arrow is adjoint to 
$$1 \sm \delta(n,k): S^{mk} \sm \C(n,k)_+ \sm S^{n(k)} \ra S^{mk} \sm S^{nk}.$$
This diagram of spectra commutes because there is a diagram of $\Sigma_k$--spaces:
\begin{equation*}
\xymatrix{
\C(m+n,k)_+ \sm (S^{m+n})^{(k)}  \ar[rr]^-{\delta(m+n,k)} & & (S^{m+n})^{\sm k} = S^{mk}\sm S^{nk} \\
\C(m+n,k)_+ \sm S^m \sm S^{n(k)} \ar[u]_{1 \sm d} & & \\
S^m \sm \C(m+n,k)_+ \sm S^{n(k)} \ar[u]_{\tau \sm 1} \ar[rr]^-{\beta(m,n,k) \sm 1}  & &
S^{mk} \sm \C(n,k)_+ \sm S^{n(k)}, \ar[uu]_{1 \sm \delta(n,k)}
}
\end{equation*}
easily checked to be commutative.
\end{proof}

\subsection{Proof of the duality theorem} \label{duality proof}

Our strategy in proving \thmref{duality thm} is to show our map is equivariantly homotopic to a duality map constructed in a standard way.

We begin with a general equivariant duality construction. Let $V$ be a real representation of a finite group $G$, and let $S^V$ denote the one point compactification $V \cup \{\infty\}$, regarded as a based $G$--space with basepoint $\infty$.

Let $\mu: V_+ \sm S^V \ra S^V$ be defined by 
\begin{equation*}
\mu(x,y) = \begin{cases} x-y & \text{if } x,y \in V \\ \infty & \text{otherwise.}
\end{cases}
\end{equation*}
This is a well--defined continuous $G$--map.

Let $K \subset S^V$ be a based $G$--subspace such that $(S^V,K)$ is an equivariant NDR pair. (Equivalently, the inclusion of $K$ into $S^V$ is an equivariant cofibration.)  The map $\mu$ induces a map of pairs:
$$ \mu: (S^V - K)_+ \sm (S^V,K) \ra (S^V, S^V - \mathbf 0).$$
The following is presumably well known.

\begin{prop} $\mu$ is an S--duality map, in the sense that 
$$ (S^V - K)_+ \sm (S^V,K) \xra{\mu} (S^V, S^V - \mathbf 0) \xla{\sim} (S^V, \infty)$$
induces a duality map $\bar{\mu}: (S^V - K)_+ \sm S^V/K \ra S^V$.
\end{prop}
\begin{proof}[Sketch Proof]  In the nonequivariant case, this can be read off of Spanier's original paper \cite{spanier}.  For the equivariant case, apply \cite[Construction III.4.5]{lmms} to the following situation: let $N \subset S^V$ be a $G$--neighborhood of $K$ such that $K \subset N$ is an equivariant homotopy equivalence, then let $X=S^V-N$ and $A = \emptyset$.  Via the equivalences $C(V,V-X) \simeq S^V/K$ and $C(X,\emptyset) \simeq (S^V-K)_+$, the map this construction yields corresponds to $\bar{\mu}$.
\end{proof}

\begin{ex} \label{duality example} This construction gives us $\Sigma_k$ -- duality maps
$$ \mu(n,k): F(\mathbb R^n, k)_+ \sm (S^{nk}, \Delta_k(S^n)) \ra (S^{nk}, S^{nk} - \mathbf 0).$$
\end{ex}

Now consider the following situation.  Suppose given
$$ \theta: \mathbb R^n_+ \sm S^n \ra S^n$$
satisfying the condition
\begin{equation} \label{cond1}
\theta(x,y) = \mathbf 0 \text{ only if } x=y.
\end{equation}
Then $\theta^{\sm k}: \mathbb R^{nk}_+ \sm S^{nk} \ra S^{nk}$ restricts to define a $\Sigma_k$--equivariant map
$$ \theta(k): F(\mathbb R^n, k)_+ \sm (S^{nk}, \Delta_k(S^n)) \ra (S^{nk}, S^{nk} - \mathbf 0).$$

\begin{ex} \label{duality example2} If  $\mu(n): \mathbb R^n_+ \sm S^n \ra S^n$ is defined by 
\begin{equation*}
\mu(n)(x,y) = \begin{cases} x-y & \text{if } x,y \in \mathbb R^n \\ \infty & \text{otherwise,}
\end{cases}
\end{equation*}
the resulting map
$$ \mu(n,k): F(\mathbb R^n, k)_+ \sm (S^{nk}, \Delta_k(S^n)) \ra (S^{nk}, S^{nk} - \mathbf 0)$$
is precisely the duality map of \exref{duality example}.
\end{ex}

A little variation on this last construction goes as follows. Suppose given
$$ \theta: \C(n,1)_+ \sm S^n \ra S^n$$
satisfying the condition
\begin{equation} \label{cond2}
\theta(c,y) = \mathbf 0 \text{ only if } c(\mathbf 0)=y.
\end{equation}
(Recall that $c(\mathbf 0)$ is the center of the little cube $c$.)
Then 
$$\theta^{\sm k}: \C(n,1)^k_+ \sm S^{nk} \ra S^{nk}$$
restricts to define a $\Sigma_k$--equivariant map
$$ \theta(k): \C(n,k)_+ \sm (S^{nk}, \Delta_k(S^n)) \ra (S^{nk}, S^{nk} - \mathbf 0).$$

\begin{ex} \label{keyex} Choose a homeomorphism $h: \stackrel{\circ}{I^n} \xra{\sim} \mathbb R^n$.  Then define $$d(n): \C(n,1)_+ \sm S^n \ra S^n$$ to be the composite
$$ \C(n,1)_+ \sm S^n \ra \stackrel{\circ}{I^n}_+ \sm S^n \xra{h \sm 1} \mathbb R^n_+ \sm S^n \xra{\mu(n)} S^n,$$
where the first map sends a little cube to its center, and $\mu(n)$ is as in \exref{duality example2}.  The resulting family of maps,
$$d(n,k): \C(n,k)_+ \sm (S^{nk}, \Delta_k(S^n)) \ra (S^{nk}, S^{nk} - \mathbf 0),$$
will be duality maps, as the maps $\mu(n,k)$ were.
\end{ex}

\begin{ex}  \label{goodex} Let $\delta(n): \C(n,1)_+ \sm S^n \ra S^n$ be defined by $\delta(n)(c,y) = c^*(y)$.  Then (\ref{cond2}) holds, and the resulting maps 
$$\delta(n,k): \C(n,k)_+ \sm (S^{nk}, \Delta_k(S^n)) \ra (S^{nk}, \infty) \hra (S^{nk}, S^{nk} - \mathbf 0)$$
are the maps of \secref{little cubes}.
\end{ex}

Let $id \in \C(n,1)$ denote the identity cube.

\begin{lem} \label{duality lem} Suppose  $\theta: \C(n,1)_+ \sm S^n \ra S^n$ satisfies (\ref{cond2}) and also
\begin{equation} \label{cond3}
 \theta(id, y) = y \text{ for all $y$ in $S^n$.} 
\end{equation}
Then $\theta(k)$ is equivariantly homotopic to the map $d(n,k)$ of \exref{keyex}.
\end{lem}

Momentarily assuming this, \thmref{duality thm} follows: since the map $\delta(n)$ of the \exref{goodex} satisfies the hypotheses of the lemma, we conclude that $\delta(n,k)$ is homotopic to the known duality map $d(n,k)$.

\begin{proof}[Proof of \lemref{duality lem}]  A map $\theta$ satisfying (\ref{cond2}) can be regarded as a map of pairs
$$ (\C(n,1)_+ \sm S^n, \C(n,1)_+ \sm S^n - \{(c,y) \ | \ c(\mathbf 0) = y\}) \xra{\theta} (S^n, S^n - \mathbf 0).$$
Now we observe that the inclusion $S^n \hra \C(n,1)_+ \sm S^n$ sending $y$ to $(id,y)$ induces a homotopy equivalence of pairs
$$ i:(S^n,S^n - \mathbf 0) \xra{\sim} (\C(n,1)_+ \sm S^n, \C(n,1)_+ \sm S^n - \{(c,y) \ | \ c(\mathbf 0) = y\}).$$
Thus maps $\theta$ satisfying (\ref{cond2}) are classified up to homotopy by the homotopy class of $\theta \circ i: (S^n, S^n - \mathbf 0) \ra (S^n, S^n - \mathbf 0)$.  But condition (\ref{cond3}) precisely says that $\theta \circ i$ is the identity map, as is $d(n) \circ i$.  Thus $\theta$ is homotopic to $d(n)$ through maps satisfying (\ref{cond2}), and so $\theta(k)$ is homotopic to $d(n,k)$ for all $k$.
\end{proof}
\section{The operad action theorem} \label{operad section}

In this section we use some of the ideas from the previous section to state and prove a more precise version of \thmref{operad theorem}.

Firstly, the structure map 
$$ \theta(r): \C(n,r) \times_{\Sigma_r} (\Omega^n X)^r \ra \Omega^nX,$$
is easy to define.  Since $\displaystyle (\Omega^n X)^r = \MapT(\bigvee_r S^n, X)$, $\theta(r)$ is the map induced by 
$$ \alpha(n,r): \C(n,r)_+ \ra \MapT(S^n, \bigvee_r S^n).$$

Next we observe that contravariant functor from spaces to towers of spectra sending $K$ to $P^K(X)$ is continuous, and that $e^K(X)$ is a natural transformation of continuous functors.  Thus one gets maps
$$\MapT(K,L) \ra \Map_{\Sp}(P^L(X),P^K(X)),$$
natural in all variables, and compatible with 
$$\MapT(K,L) \ra \MapT(\MapT(L,X), \MapT(K,X)).$$
Adjointing the $\Sigma_r$--equivariant composite
$$ \C(n,r)_+ \xra{\alpha(n,r)} \MapT(S^n, \bigvee_r S^n) \ra \Map_{\Sp}(P^{\bigvee_r S^n}(X),P^{S^n}(X))$$
yields a natural maps of towers
$$ \theta(r): \C(n,r)_+ \sm _{\Sigma_r} P^{\bigvee_r S^n}(X) \ra P^{S^n}(X).$$
This is the map of \thmref{operad theorem}, and property (1) listed there clearly holds.

It remains to identify the map on $k^{th}$ fibers in terms of the little $n$--cubes operad structure.  

One approach to the operad structure is as follows. Let $\C(n,r,k)$ be the space of all embeddings
$$ c: \coprod_{i=1}^k \stackrel{\circ}{I^n} \ra \coprod_{j=1}^r \stackrel{\circ}{I^n}$$
such that each nontrivial component is a little $n$--cube.  This is a $\Sigma_r \times \Sigma_k$--space.  Then the operad structure is given by the $\Sigma_k$--equivariant maps
\begin{equation*}  \circ_{r,k}: \C(n,r) \times_{\Sigma_r} \C(n,r,k) \ra \C(n,k)
\end{equation*}
sending $(c,d)$ to the composition $c\circ d$.

To see that this agrees with the definition given in \cite{may}, note that $\C(n,r,k)$ decomposes as the product, over all maps $\lambda: \mathbf k \ra \mathbf r$, of the spaces
$$ \C(n,k_1) \times \dots \times \C(n,k_r),$$
where $k_j = |\lambda^{-1}(j)|$, and the corresponding component of $\circ_{r,k}$ is the usual structure map.
 
The $k^{th}$ fiber of $P^{\bigvee_r S^n}(X)$ is naturally equivalent to $$\MapS((\bigvee_r S^n)^{(k)}, X^{\sm k})_{h\Sigma_k},$$ and, by construction, the map $\theta(r,k)$ induced by $\theta(r)$ on $k^{th}$ fibers is induced in the apparent way by $\alpha(n,r)$.  

Generalizing definitions in \secref{duality section}, let
$$\delta(n,r,1): \C(n,r,1)_+ \sm (\bigvee_r S^n) \ra S^n$$
be the map sending $(c,t)$ to $c^*(t)$, where $c^*: \bigvee_r S^n \ra S^n$ is the Thom--Pontryagin collapse map associated to $c$.  As before, $\delta(n,r,1)^{\sm k}$ induces maps
$$ \delta(n,r,k): \C(n,r,k)_+ \sm (\bigvee_r S^n)^{(k)} \ra S^{nk}.$$
Using \thmref{duality thm}, one can deduce that the stable adjoint of $\delta(n,r,k)$, 
$$ \Tilde{\delta}(n,r,k): \Sinfty \C(n,r,k)_+ \ra \MapS((\bigvee_rS^n)^{(k)}, S^{nk})$$
is an equivalence.  There results a natural weak equivalence of $\Sigma_r$--spectra:
$$ \C(n,r,k)_+ \sm_{\Sigma_k} (\Sigma^{-n} X)^{\sm k} \simeq \MapS((\bigvee_r S^n)^{(k)}, X^{\sm k})_{h\Sigma_k}.$$
We note that it is easy to see that this equivalence is compatible with the decomposition of $\MapS((\bigvee_r S^n)^{(k)}, X^{\sm k})_{h\Sigma_k}$ arising in our product theorems.

Property (2) of \thmref{operad theorem} can now be more precisely stated.

\begin{prop}  There is a commutative diagram of weak natural transformations:
{
\begin{equation*}
\xymatrix{
\C(n,r)_+ \sm_{\Sigma_r} \C(n,r,k)_+ \sm_{\Sigma_k} (\Sigma^{-n} X)^{\sm k}  \ar[rr]^-{\circ_{r,k} \sm 1} \ar[d]^-{\wr} &&
\C(n,k)_+ \sm_{\Sigma_k} (\Sigma^{-n} X)^{\sm k} \ar[d]^-{\wr}  \\
 \C(n,r)_+ \sm_{\Sigma_r} \MapS((\bigvee_r S^n)^{(k)}, X^{\sm k})_{h\Sigma_k} \ar[rr]^-{\theta(r,k)(X)}&&
\MapS(S^{n(k)}, X^{\sm k})_{h\Sigma_k}.
}
\end{equation*}}
\end{prop}

\proof It suffices to show that there is a commutative diagram of $\Sigma_k$--spaces
\begin{equation*}
\xymatrix{
\C(n,r) \times_{\Sigma_r} \C(n,r,k)  \ar[rr]^-{\circ_{r,k}} \ar[d]^-{1 \times \Tilde{\delta}(n,r,k)} && \C(n,k)  \ar[d]^-{\Tilde{\delta}(n,k)}
  \\
 \C(n,r) \times_{\Sigma_r} \MapT((\bigvee_r S^n)^{(k)}, S^{nk}) \ar[rr]^-{\alpha} &&
\MapT(S^{n(k)}, S^{nk})
}
\end{equation*}
where $\alpha$ is induced by $\alpha(n,r)$.
That this diagram commutes is easily checked: given $c \in \C(n,r)$ and $d = (d_1, \dots, d_k) \in \C(n,r,k)$, we have that  
$$(\alpha \circ (1 \times \Tilde{\delta}(n,r,k)))(c,d) = (d_1^* \circ c^*) \sm \dots \sm (d_k^* \circ c^*), $$
while 
$$(\Tilde{\delta}(n,k)\circ  (\circ_{r,k}))(c,d) = (c \circ d_1)^*\sm \dots \sm (c \circ d_k)^*.$$
These agree because the Thom--Pontryagin collapse is a contravariant functor: $$(c\circ d)^* = d^* \circ c^*.\eqno{\qed}$$

\appendix
\section{A proof of Arone's theorem when $K$ is a sphere} \label{convergence appendix}

In this appendix, we give a proof of \thmref{arone thm} in the case when $K$ is a sphere.  The first reductions follow along the lines of Arone's proof in \cite{arone}, but we use \thmref{duality thm} to simplify the proof of the last key step: proving that a certain `cross effect' map is an equivalence.  It seems likely that some variant of our proof can be used to prove the theorem in general.  As does Arone, we use ideas from Goodwillie calculus at a couple of points.  

The theorem we are trying to prove says that 
$$ e_{\infty}^{S^n}(X): \Sinfty \MapT(S^n,X)_+ \ra P_{\infty}^{S^n}(X)$$
is a weak homotopy equivalence if the connectivity of $X$ is greater than $n$.

To explain why this is equivalent to \thmref{arone thm}, and to effect our first reduction, let $F(X)$ and $G(X)$ respectively denote the domain and range of $ e_{\infty}^{S^n}(X)$.  \cite[Example 4.5]{goodwillie2} says that $F(X)$ is $n$--analytic: this means that $F$ behaves well (in a precise sense defined in \cite{goodwillie2}) on $n$-connected spaces.  Consideration of the fibers of the Arone tower, identified in \propref{fiber prop}, shows that the projection maps
$$ q_k: G(X) = P_{\infty}^{S^n}(X) \ra P_{k}^{S^n}(X)$$
are  $(1+ \text{conn }X - n)(1+k)-1$ connected, and then that  $G(X)$ is $n$-analytic. 

Goodwillie's argument proving  \cite[Prop.5.1]{goodwillie2} then shows that $e_{\infty}^{S^n}(X)$ will be a weak equivalence for all $n$--connected $X$ if it is an equivalence for all $X$ of the form $\Sigma^n Y$, with $Y$ connected.

At this point, we use the classical model $C_n(Y)$ for $\MapT(S^n, \Sigma^n Y)$ built from the spaces $\C(n,k)$.  

For a space $Y$ with basepoint $*$, let
$$ C_n(Y) = (\coprod_{k = 1}^{\infty} \C(n,k)\times_{\Sigma_k} Y^k)/({\sim}),$$
where $(c_1,\dots, c_k, y_1, \dots, y_{k-1}, *) \sim (c_1,\dots, c_{k-1}, y_1, \dots, y_{k-1})$ generates the equivalence relation.  $C_n(Y)$ is filtered, with 
$$ F_kC_n(Y) = (\coprod_{j = 1}^{k} \C(n,j)\times_{\Sigma_j} Y^j)/({\sim}),$$
and there are cofibration sequences
\begin{equation*} 
F_{k-1}C_n(Y) \ra F_kC_n(Y) \ra \C(n,k)_+ \sm_{\Sigma_k} Y^{\sm k}.
\end{equation*}
A natural map \ $ \alpha_n(Y): C_n(Y) \ra \MapT(S^n, \Sigma^n Y)$
is defined to be the map induced by the composites
{\small $$\C(n,k) \times Y^k \xra{\alpha(n,k)\times \eta} \MapT(S^n, \bigvee_k S^n) \times \MapT(\bigvee_k S^n, \Sigma^n Y) \xra{\circ} \MapT(S^n, \Sigma^n Y).$$}
Explicitly, if we regard $c \in \C(n,k)$ as an embedding $c: \coprod_{k} \stackrel{\circ}{I^n} \longrightarrow \stackrel{\circ}{I^n}$, and $y \in Y^k$ as a based map $y: \bigvee_k S^0 \ra Y$, then 
$$ \alpha_n(Y)([c,y]) = (\Sigma^n y) \circ c^*.$$
The classical theorem, \cite[Thm.2.7]{may}, then states that $\alpha_n(Y)$ is a weak homotopy equivalence if $Y$ is connected.  

It follows that, to show $e_{\infty}^{S^n}(\Sigma^n Y)$ is an equivalence, it suffices to show that 
$$\theta(n,k,Y): \Sinfty F_kC_n(Y)_+ \ra P^{S^n}_k(\Sigma^n Y)$$
is an equivalence, where $\theta(n,k,Y)$ is the composite
\begin{equation*}
\begin{split} \Sinfty F_kC_n(Y)_+  \hra \Sinfty C_n(Y)_+ & \xra{\alpha_n(Y)} \Sinfty \MapT(S^n, \Sigma^n Y)_+ \\ & \xra{e_k^{S^n}(\Sigma^n Y)} P^{S^n}_k(\Sigma^n Y).
\end{split}
\end{equation*} 
This we proceed to show by induction on $k$, using ideas from \cite{goodwillie3}.  If $G: \T \ra \Sp$ is a functor, we let $p_kG: \T \ra \Sp$ be its universal $k$--excisive quotient, and $\chi_kG: \T^k \ra \Sp$ its $k^{th}$ cross effect.

It is quite easy to see that $p_k$ preserves fibration sequences of functors, and, since our functors take values in spectra, also cofibration sequences.  Also, the functor $\C(n,k)_+ \sm_{\Sigma_k} Y^{\sm k}$ is an example of a homogeneous functor of degree $k$, i.e.\ a functor $G$ with $G \simeq p_kG$ and $p_{k-1}G \simeq *$. 

Such considerations show that both $\Sinfty F_kC_n$ and $P^{S^n}_k$ are $k$--excisive, and both $i:F_{k-1}C_n \hra F_kC_n$ and $p:P^{S^n}_k \ra P^{S^n}_{k-1}$ induce equivalences after applying $p_{k-1}$.  Then the inductive hypothesis, combined with the  commutativity of the diagram
\begin{equation*}
\xymatrix{
\Sinfty F_kC_n(Y)_+ \ar[rr]^-{\theta(n,k, Y)} &&
P^{S^n}_k(\Sigma^nY) \ar[d]^p  \\
\Sinfty F_{k-1}C_n(Y)_+  \ar[rr]^-{\theta(n,k-1,Y)} \ar[u]_i &&
P^{S^n}_{k-1}(\Sigma^nY)
}
\end{equation*}
shows that $p_{k-1}\theta(n,k,Y)$ is an equivalence for all $Y$.

\cite[Proposition 3.4]{goodwillie3} implies that a natural map $\theta$ between $k$--excisive functors will be an equivalence if both $p_{k-1}\theta$ and $\chi_k\theta$ are equivalences.  Thus, our inductive proof that $\theta(n,k,Y)$ is an equivalence will be complete if we establish that 
$$ \chi_k \theta(n,k)(Y_1,\dots,Y_k):\chi_k \Sinfty F_kC_n(Y_1,\dots,Y_k) \ra
\chi_kP^{S^n}_k(\Sigma^nY_1,\dots, \Sigma^nY_k)$$
is an equivalence.

Cross effects are defined as iterated fibers of certain cubical diagrams.  But if a functor $G$ takes values in spectra (where finite coproducts are equivalent to finite products), $\chi_kG(Y_1,\dots,Y_k)$ is naturally equivalent to the iterated {\em cofiber} of the $k$--dimensional cube with entries $G(Z_1\vee \dots \vee Z_k)$ with $Z_j \in \{*,Y_j\}$, and the canonical map $G(Y_1 \vee \dots \vee Y_k) \ra \chi_kG(Y_1, \dots, Y_k)$ is a retraction.

It is immediate from the definitions that, when $G = \Sinfty F_kC_n$, this  retraction is the map
$$ \Sinfty F_kC_n(Y_1 \vee \dots \vee Y_k)_+ \ra \Sinfty \C(n,k)_+ \sm (Y_1 \sm \dots Y_k)$$
induced by the evident map
$$\C(n,k)_+ \times_{\Sigma_k} (Y_1 \vee \dots \vee Y_k)^k \ra \C(n,k)_+ \sm (Y_1 \sm \dots \sm Y_k).$$
Now note that the functor $\rho: \Sp^{\E_k} \ra \Sp^{\Sigma_k}$, which sends a functor $X$ to the $\Sigma_k$--spectrum $X_k/X_{k-1}(\mathbf k)$, induces a natural map 
$$ \rho: \Map_{\Sp}^{\E_k}(K^{\sm}, X^{\sm}) \ra \Map_{\Sp}^{\Sigma_k}(K^{(k)}, X^{(k)})$$
such that 
$$\Map^{\Sigma_n}_{\T}(K^{(k)}, X^{\sm k}) \xra{i}\Map_{\Sp}^{\E_k}(K^{\sm}, X^{\sm}) \xra{\rho} \Map{\Sp}^{\Sigma_k}(K^{(k)}, X^{(k)})$$
is induced by $\pi:X^{\sm k} \ra X^{(k)}$.  It follows that, when $G = P^{S^n}_k \circ \Sigma^n$, the retraction is equivalent to the composite
\begin{equation*}
\begin{split} \Map_{\Sp}^{\E_k}(S^{\sm},(\Sigma^n Y_1 \vee \dots \vee \Sigma^n Y_k)^{\sm}) & \xra{\rho} \Map_{\Sp}^{\Sigma_k}(S^{n(k)},(\Sigma^n Y_1 \vee \dots \vee \Sigma^n Y_k)^{(k)}) \\
& \ra \MapS(S^{n(k)}, \Sigma^n Y_1 \sm \dots \Sigma^n Y_k).
\end{split}
\end{equation*}

If $Y$ is a space, let $\Tilde{\delta}(n,k,Y): \C(n,k)_+ \sm Y \ra \MapS(S^{n(k)}, S^{nk} \sm Y)$ be the map adjoint to $\delta(n,k) \sm 1_Y$. By \thmref{duality thm}, $\Tilde{\delta}(n,k,Y)$ will be a weak equivalence.

\begin{lem} The diagram
{\footnotesize 
\begin{equation*}
\xymatrix{
\Sinfty \C(n,k)_+ \times_{\Sigma_k} (Y_1 \vee \dots \vee Y_k)^k \ar[d] \ar[rrr]^-{\theta(n,k)(Y_1\vee \dots \vee Y_k)} &&&
\Map_{\Sp}^{\E_k}(S^{\sm},(\Sigma^n Y_1 \vee \dots \vee \Sigma^n Y_k)^{\sm}) \ar[d]  \\
\Sinfty \C(n,k)_+ \sm (Y_1 \sm \dots \sm Y_k) \ar[rrr]^{\Tilde{\delta}(n,k,Y_1 \sm \dots \sm  Y_k)} &&&
\MapS(S^{n(k)}, \Sigma^n Y_1 \sm \dots \sm \Sigma^n Y_k)
}
\end{equation*}}
commutes.

\end{lem}

Assuming this, our proof of \thmref{arone thm} is done, as we can identify the map $\chi_k \theta(n,k)(Y_1,\dots,Y_k)$ with the weak equivalence $\Tilde{\delta}(n,k,Y_1 \sm \dots \sm Y_k)$.

The functor $\T \ra \T$ sending $Y$ to $Y^{(k)}$ is continuous.  Thus it induces natural maps $\xi^{(k)}: \MapT(K,Y) \ra \MapT(K^{(k)}, Y^{(k)})$.  Recalling the definition of the maps $e_k^K(X)$, one sees that the commutativity of the diagram in the lemma will follow from the commutativity of the following diagram of spaces, when $Y = Y_1 \vee \dots \vee Y_k$:
{
\begin{equation*}
\xymatrix{
\C(n,k)_+ \times Y^k \ar[ddd] \ar[rr]^-{\alpha(n,k) \times \eta} &&  \MapT(S^n,\bigvee_k S^n) \times \MapT(\bigvee_k S^n, \Sigma^n Y) \ar[d]^{\circ}  \\
&& \MapT(S^n, \Sigma^n Y) \ar[d]^{\xi^{(k)}} \\
&& \MapT(S^{n(k)}, (\Sigma^n Y)^{(k)}) \ar[d] \\
\C(n,k)_+ \sm (Y_1 \sm \dots \sm Y_k) \ar[rr]^-{\Tilde{\delta}(n,k,Y_1 \sm \dots \sm Y_k)} && \MapT(S^{n(k)}, \Sigma^n Y_1 \sm \dots \Sigma^n Y_k).
}
\end{equation*}}
The commutativity of this diagram is verified easily.  That the two maps agree on an element $(c_1, \dots, c_k ,y_1,\dots,y_k)$, with $c_i\in \C(n,1)$ and $y_i \in Y_i$ (viewed as a map $y_i: S^0 \ra Y_i$), amounts to the observation that the diagram
\begin{equation*}
\xymatrix{
S^{nk} \ar[d]^{c_1^* \sm \dots \sm c_k^*} \ar[r]^-{c^{\sm k}} &
(\bigvee_k S^n)^{\sm k} \ar[rrr]^-{(\Sigma^n(y_1 \vee \dots \vee y_k))^{\sm k}} &&& (\Sigma^n Y_1 \vee \dots \vee \Sigma^n Y_k)^{\sm k} \ar[d]  \\
S^{nk} \ar[rrrr]^-{\Sigma^n y_1 \sm \dots \sm \Sigma^n y_k} &&&&
\Sigma^n Y_1 \sm \dots \sm \Sigma^n Y_k
}
\end{equation*}
commutes.

\begin{rem}  The argument here shows that the natural transformations
$$ e^{S^n}_k(\Sigma^n Y): \Sinfty (\Omega^n \Sigma^n Y)_+ \ra P^{S^n}_k(\Sigma^n Y)$$
and 
$$ p^{S^n}_k(\Sigma^n Y): P^{S^n}_k(\Sigma^n Y)\ra P^{S^n}_{k-1}(\Sigma^n Y)$$
admit compatible natural weak right inverses
$$ s_k: P^{S^n}_k(\Sigma^n Y) \ra \Sinfty (\Omega^n \Sigma^n Y)_+ $$
and 
$$ t_k:  P^{S^n}_{k-1}(\Sigma^n Y)\ra P^{S^n}_{k}(\Sigma^n Y).$$
This is a form of the Snaith splitting theorem, and the splitting obtained this way is equivalent to other classical constructions \cite[Appendix B]{k2}.  \\

\end{rem}

\Addresses\recd
\end{document}